\RequirePackage{fix-cm}
\documentclass[smallextended]{svjour3}       
\smartqed  
\usepackage{color}
\usepackage{bm}
\usepackage{amsmath}
\usepackage{latexsym,amssymb,amsxtra,amsfonts}
\usepackage[cmtip,all]{xy}
\usepackage{amsfonts}
\usepackage{amssymb}
\usepackage{url}
\usepackage{bm}
\usepackage{lineno} 
\usepackage[footnotesize,FIGBOTCAP,TABTOPCAP]{subfigure}
\usepackage[font=footnotesize, labelfont=bf]{caption} 
\usepackage{graphics}
\usepackage{comment}
\usepackage{graphicx}
\usepackage{multirow}
\usepackage{chngcntr}
\usepackage{setspace}
\usepackage{marginnote}
\usepackage{threeparttable}
\usepackage{enumerate}
\usepackage[colorlinks=true, bookmarksopen,
            pdfauthor={Hoang Tran},
            pdfcreator={pdftex},
            pdfsubject={algorithms},
            linkcolor={blue},
            anchorcolor={black},
            citecolor={red},
            filecolor={magenta},
            menucolor={black},
            pagecolor={red},
            plainpages=false,pdfpagelabels,
            urlcolor={blue}]{hyperref}
\usepackage{cite}

\newtheorem{assumption}{Assumption}

\newcommand{\sproof}{{\em Proof.}  }
\newcommand{\fproof}{\hfill $\square$}
\numberwithin{equation}{section}

\graphicspath{{./}}

\newcommand{\Loeve}{Lo{\`e}ve }

\begin{document}

\title{Analysis of quasi-optimal polynomial approximations for parameterized {PDEs} with deterministic and stochastic coefficients
}
\subtitle{}

\titlerunning{Analysis of quasi-optimal approximations for parameterized {PDEs}}        

\author{
Hoang Tran \and 
Clayton G.~Webster \and 
Guannan Zhang
}


\institute{
           Hoang Tran \at
              Department of Computational and Applied Mathematics, Oak Ridge National Laboratory
              \\
               1 Bethel Valley Road, P.O. Box 2008, Oak Ridge TN 37831-6164 
              \\
               \email{tranha@ornl.gov}
           \and 
            Clayton G.~Webster \at
              Department of Computational and Applied Mathematics, Oak Ridge National Laboratory
              \\
               1 Bethel Valley Road, P.O. Box 2008, Oak Ridge TN 37831-6164 
              \\
               \email{webstercg@ornl.gov}
           \and 
          Guannan Zhang \at
             Department of Computational and Applied Mathematics, Oak Ridge National Laboratory
              \\
              1 Bethel Valley Road, P.O. Box 2008, Oak Ridge TN 37831-6164 
              \\
               \email{zhangg@ornl.gov}    
  }

\date{}

\maketitle

\begin{abstract}
In this work, we present a generalized methodology for analyzing the convergence of quasi-optimal Taylor and Legendre approximations, applicable to a wide class of parameterized elliptic PDEs with {  finite-dimensional} deterministic and stochastic inputs.  
Such methods construct an {  optimal index set that corresponds to the sharp estimates of the polynomial coefficients.}
%
%
{  Our analysis furthermore represents a new approach for estimating best $M$-term approximation errors by means of coefficient bounds, without using Stechkin inequality.} 
%
%
The framework we propose for analyzing asymptotic truncation errors is based on 
an extension of the underlying multi-index set into a continuous domain, and then an 
approximation of the cardinality (number of integer multi-indices) by its Lebesgue measure.
%
%
%
%
Several types of isotropic and anisotropic (weighted) multi-index sets are explored, and  
rigorous proofs reveal sharp asymptotic error estimates in which we achieve sub-exponential convergence rates 
(of the form 
$M \textrm{exp}({-(\kappa M)^{1/N}})$, with $\kappa$ a constant depending on the shape and size of multi-index sets) 
with respect to the total number of degrees of freedom.
%
Through several theoretical examples, we explicitly derive the constant $\kappa$  
and use the resulting sharp bounds to illustrate the effectiveness of Legendre 
over Taylor approximations, as well as compare our rates of convergence with current published results. Computational evidence complements the theory and shows
the advantage of our generalized framework compared to previously developed estimates. 
\end{abstract}

\section{Introduction}
\label{sec:intro}
This paper focuses on a relevant model boundary value problem, involving the simultaneous solution of 
a family of equations, parameterized by a {  finite-dimensional} vector 
${\bm y}=(y_1, \ldots, y_N)\in\Gamma=\prod_{i=1}^N\Gamma_i\subset\mathbb{R}^N$, 
on a bounded Lipschitz domain $D\subset\mathbb{R}^d,\, d\in \{1,\,2,\,3\}$.  In particular, we consider 
a differential operator $\mathcal{L}$ defined on $D$, and let $a(x, {\bm{y}})$, with $x\in D$ and 
${\bm y}\in\Gamma$, represent the input coefficient associated with the operator $\mathcal{L}$.
The forcing term $f=f(x)\in L^2(D)$ is assumed to be a fixed function of $x\in D$.
We concentrate on the following parameterized boundary value problem:  
for all ${\bm y}\in\Gamma$, find $u(\cdot, {\bm y}):{\overline{D}}\rightarrow\mathbb{R}$, 
such that the following equation holds
\begin{equation}
\label{eq:genPDE}
\mathcal{L}(a(\cdot, {\bm{y}}))
\left[u(\cdot,{\bm{y}})\right] = f(\cdot) \quad\mbox{in } D,
\end{equation}
subject to suitable (possibly parameterized) boundary conditions.   
%
%
We require $a$ and $f$ to be chosen such that system \eqref{eq:genPDE} is well-posed in a Banach space,
with unique solution $u$, 
such that, when suppressing the explicit dependence on $x$, the map ${\bm y}\mapsto u({\bm y})$ is defined from the parameter domain $\Gamma$ into the solution space $V(D)$.

Problems such as \eqref{eq:genPDE} arise in contexts of both deterministic and stochastic modeling. 
In the deterministic setting, the parameter vector 
${\bm y}$ is known or controlled by the user, and a typical goal is to study the dependence of $u$ on these parameters, e.g., optimizing an output of the equation with respect to ${\bm y}$ 
(see \cite{Buffa:2012iz,Milani:2008hm} for more details).  
On the other hand, stochastic modeling is motivated by many engineering and science problems in which the input data is not known exactly. A quantification of the effect of the input uncertainties on the output of simulations is necessary to obtain a reliable prediction of the physical system.
A natural way to incorporate the presence of input uncertainties into the governing model \eqref{eq:genPDE} is to consider the parameters $\{y_n(\omega)\}_{n=1}^N$ as random variables and ${\bm y}(\omega):\Omega\rightarrow\Gamma$ a 
random vector, where $\omega\in\Omega$ and $\Omega$ is the set of outcomes. In this setting, we assume the components 
of ${\bm y}$ have a joint probability density function (PDF) $\varrho:\Gamma\rightarrow \mathbb{R}_+$, with 
$\varrho\in L^{\infty}(\Gamma)$ known directly through, e.g., 
truncations of correlated random fields \cite{MR0651017,MR0651018,Wiener_38,Ghanem_Spanos_91},
%
such that the probability space is equivalent to $(\Gamma, \mathcal{B}(\Gamma), \varrho(\bm{y})d\bm{y})$, where $\mathcal{B}(\Gamma)$ denotes the Borel $\sigma$-algebra on $\Gamma$ and $\varrho(\bm{y})d\bm{y}$ is the probability measure of $\bm{y}$.  

Monte Carlo (MC) methods (see, e.g., \cite{Fis96}) are the most popular approaches for approximating high-dimensional 
integrals, such as expectation or two-point correlation, based on independent realizations $u({\bm y}_k)$, $k=1,\ldots, M$, of the solution to \eqref{eq:genPDE}; 
approximations of the expectation or other QoIs are obtained by averaging over the corresponding realizations of that quantity. The resulting numerical error is proportional to $M^{-1/2}$, 
thus, achieving convergence rates independent of dimension $N$, but requiring a very large number of samples to achieve reasonably small errors.    
Moreover, MC methods do not have the ability to simultaneously approximate the solution map ${\bm y}\mapsto u({\bm y})$, since they are quadrature techniques and do not exploit the fact in many scenarios, the solutions smoothly depend on the coefficient $a$
Taking this smooth dependence into account, several global polynomial approximation techniques, for instance, {\em intrusive} Galerkin methods \cite{BTZ04,TS07} and {\em non-intrusive} collocation methods \cite{BNT07,NTW08}, have been proposed, often featuring much faster convergence rates.

Let $\mathcal{S} = \left\{{\bm {\nu}} = (\nu_i)_{1\le i\le N}:\nu_i \in\mathbb{N} \right\}$. 
Global polynomial approximation methods seek to build an approximation $u_\Lambda$ to the solution $u$ of the form:
 \begin{align}
 \label{global_poly}
 u_\Lambda (x,{\bm y}) = \sum_{{\bm \nu}\in\Lambda} c_{\bm \nu}(x) {{\bm \Psi}}_{\bm \nu} ({\bm y}),
 \end{align}  
for a finite multi-index set $\Lambda\subset \mathcal{S}$, where ${{\bm \Psi}}_{\bm \nu}$ is a multivariate polynomial in ${span}\{{\bm y}^{\bm \mu}:{\bm \mu} \le {\bm \nu}\}$ for ${\bm \nu} \in \Lambda$ and $c_{\bm \nu} \in V(D)$ is the coefficient to be computed, both of which are method specific. Here, for two vectors ${\bm \nu},{\bm \mu}\in\mathcal{S}$, we say ${\bm \mu} \le {\bm \nu}$ if and only if $\mu_i\le \nu_i$ for all $1\le i\le N$. Also, given ${\bm \alpha}=(\alpha_i)_{1\le i\le N}$ a vector of real numbers, we define 
$ {\bm \alpha}^{\bm \nu} = \prod_{1\le i\le N} \alpha_i^{\nu_i}$ with the convention $0^0 :=1$. We will often suppress the dependence on $x$ and use the notations $u({\bm y}):=u(\cdot,{\bm y})$ and $a({\bm y}):=a(\cdot,{\bm y})$ without loss of generality. In this paper, we are interested in solving \eqref{eq:genPDE} using a class of polynomial approximations based on the Taylor and Legendre expansions of solution $u$. The polynomial basis considered herein is thus given by the monomials ${\bm \Psi}_{\bm \nu}({\bm y}) = {\bm y}^{\bm \nu}$ (in the former case) and the tensorized Legendre polynomials ${\bm \Psi}_{\bm \nu}({\bm y}) = L_{\bm \nu}({\bm y})$ (in the latter case).



The evaluation of $u_\Lambda$ in \eqref{global_poly} requires the computation of $\# (\Lambda)$ coefficients $c_{\bm \nu}(x) \in V(D)$, where $\# (\Lambda)$ is the cardinality of $\Lambda$. A naive choice of $\Lambda$ and their corresponding polynomial spaces $\mathbb{P}_{\Lambda}(\Gamma) = span \{{\bm \Psi}_{\bm \nu}({\bm y}),\, {\bm \nu} \in \Lambda\}$, for instance, tensor product polynomial spaces, could lead to an infeasible computational cost, especially when the dimension of the parameter domain is high. It is important to be able to construct the set of the most effective indices for the approximation \eqref{global_poly}, which provides maximum accuracy for a given cardinality. In other words, given a fixed $M\in\mathbb{N}$, one searches for a set {  $\Lambda$ which minimizes the error $u - \sum_{{\bm \nu}\in\Lambda} c_{\bm \nu} {{\bm \Psi}}_{\bm \nu}$ among all index subsets of $\mathcal{S} $ of cardinality $M$}. This practice has been known as \textit{best M-term approximations}.  

The literature on the best $M$-term Taylor and Galerkin approximations has been growing fast recently, among them we refer to \cite{BTNT12,BNTT14,BAS09,CCDS13,CDS10,CDS11,CCS14,HS13,HS13b,HS13c}. In the benchmark work \cite{CDS11}, the analytic dependence of the solutions of parametric elliptic PDEs on the parameters was proved under mild assumptions on the input coefficients, and convergence analysis of the best $M$-term Taylor and Legendre approximations was established subsequently. Consider, for example, the expansion of $u$ on $\Gamma = [-1,1]^{{N}}$ by a family of $L^{\infty}$ normalized polynomials, i.e., $\|\bm\Psi_{\bm \nu}\|_{L^{\infty}(\Gamma)} = 1$. Application of the triangle inequality yields
\begin{align}
\label{triangle1}
\sup_{{\bm y}\in \Gamma}\left\|u( {\bm y}) - \sum_{{\bm \nu}\in\Lambda} c_{\bm \nu}{{\bm \Psi}}_{\bm \nu}({\bm y})\right\|_{V(D)} \le \sum_{{\bm \nu}\in \Lambda^c}\|c_{\bm \nu}\|_{V(D)}, 
\end{align}
which suggests determining the optimal index set {  $\Lambda^{\textbf{best}}_M$} by choosing the set of indices ${\bm \nu}$ corresponding to $M$ largest $\|c_{\bm \nu}\|_{V(D)}$. Here, $\Lambda^c$ denotes the complement of $\Lambda$ in $\mathcal{S}$. In \cite{CDS11}, the error of such approximation was estimated due to Stechkin inequality (see, e.g., \cite{DeV98}) such that
\begin{align}
\label{stechkin}
 \sum_{{\bm \nu}\in (\Lambda^{\textbf{best}}_M)^c}\|c_{\bm \nu}\|_{V(D)}  \le \|(\|c_{\bm \nu}\|_{V(D)})\|_{\ell^p(\mathcal{S})} M^{1-\frac{1}{p}}, 
\end{align}
where $p$ is some number in $(0,1)$ such that $(\|c_{\bm \nu}\|_{V(D)})_{\bm \nu \in \mathcal{S}}$ is $\ell^p$-summable. It should be noted that the convergence rate \eqref{stechkin} does not depend on the dimension of the parameter domain $\Gamma$ (which is possibly countably infinite therein). This error estimate, however, has some limitations. First, explicit evaluation of the coefficient $\|(\|c_{\bm \nu}\|_{V(D)})\|_{\ell^p(\mathcal{S})}$ is inaccessible in general (thus so is the total estimate). Secondly, \eqref{stechkin} often occurs with infinitely many values of $p$ and stronger rates, corresponding to smaller $p$, are also attached to bigger coefficients. For a specific range of $M$, the effective rate of convergence is unclear. 
In implementation, finding the best index set and polynomial space is an infeasible task, since this requires computation of all of the $c_{\bm \nu}$. As a strategy to circumvent this challenge, adaptive algorithms which generate the index set in a near optimal, greedy procedure were developed in \cite{CCDS13}. This method {  however comes} with a high cost of exploring the polynomial space, which may be daunting in high-dimensional problems.

Instead of building the index set based on exact values of polynomial coefficients $c_{\bm \nu}$, an attractive alternative approach (referred to as \textit{quasi-optimal approximation} throughout this paper) is to establish sharp upper bounds of $c_{\bm \nu}$ (by a priori or a posteriori methods), and then construct {  the index set $\Lambda^{\textbf{qsi}}_M$} corresponding to $M$ largest such bounds. For this strategy, the main computational work for the selection of the (near) best terms reduces to determining sharp coefficient estimates, which is expected to be significantly cheaper than exact calculations. Quasi-optimal polynomial approximation has been performed for some parametric elliptic models with optimistic results: while the upper bounds of $\|c_{\bm \nu}\|_{V(D)}$ (denoted from now by $B({\bm \nu})$) were computed with a negligible cost, the method was comparably as accurate as best $M$-term approach, as shown in \cite{BTNT12,BNTT14}. The first rigorous numerical analysis of quasi-optimal approximation was presented in \cite{BNTT14} for $B({\bm \nu}) = {\bm \rho}^{- {\bm \nu}}$ with ${\bm \rho} $ being a vector $(\rho_i)_{1\le i\le N}$ with $\rho_i > 1\ \, \forall i$. In that work, the asymptotic sub-exponential convergence rate was proved based on optimizing the Stechkin estimation. Briefly, the analysis applied Stechkin inequality to yield 
\begin{align}
\label{stechkin2}
\sum_{{\bm \nu}\in (\Lambda^{\textbf{qsi}}_M)^c}B({\bm \nu}) \le \|B({\bm \nu})\|_{\ell^p(\mathcal{S})} M^{1-\frac{1}{p}},
\end{align}
then took advantage of the formula of $B({\bm \nu})$ to compute $p\in (0,1)$, depending on $M$, which minimizes $\|B({\bm \nu})\|_{\ell^p(\mathcal{S})} M^{1-\frac{1}{p}}$. 

Although known as an essential tool to study the convergence rate of best $M$-term approximations, Stechkin inequality is probably less efficient for quasi-optimal methods. As a generic estimate, it does not fully exploit the available information of the decay of coefficient bounds. In such a setting, a direct estimate of $\sum_{{\bm \nu}\in (\Lambda^{\textbf{qsi}}_M)^c}B({\bm \nu}) $ may be viable and advantageous to provide an improved result. In addition, the process of solving the minimization problem $p^* = \mathrm{argmin}_{ {p\in(0,1)} } \|B({\bm \nu})\|_{\ell^p(\mathcal{S})} M^{1-\frac{1}{p}}$ needs to be tailored to $B({\bm \nu})$, making this approach not ideal for generalization. Currently, this minimization approach has been limited for some quite simple types of upper bounds. In many scenarios, the sharp estimates of the coefficients may involve complicated bounds which are not even explicitly computable, such as those proposed in \cite{CDS11}. The extension of this approach to such cases seems to be impossible.  

In this work, we present a generalized methodology for convergence analysis of quasi-optimal polynomial approximations for parameterized PDEs with deterministic and stochastic coefficients. {  As the errors of best $M$-term approximations are bounded by those of quasi-optimal methods
\begin{align}
\label{accessible}
\sum_{{\bm \nu}\in (\Lambda^{\textbf{best}}_M)^c}\|c_{\bm \nu}\|_{V(D)}  \le \sum_{{\bm \nu}\in (\Lambda^{\textbf{qsi}}_M)^c}\|c_{\bm \nu}\|_{V(D)} \le \sum_{{\bm \nu}\in (\Lambda^{\textbf{qsi}}_M)^c}B({\bm \nu}), 
\end{align}
our analysis also gives \textit{accessible} estimates (i.e., estimates depending only on {known} or computable terms) for best $M$-term approximation errors under several established properties on the decaying of the polynomial coefficients. These are sharp explicit theoretical estimates in the scenario that: 1) the triangle inequality \eqref{triangle1} has to be employed; 2) one has to evaluate $c_{\bm \nu}$ via their bounds $B({\bm \nu})$, in which case \eqref{accessible} represents the smallest accessible bound of $\sum_{{\bm \nu}\in (\Lambda^{\textbf{best}}_M)^c}\|c_{\bm \nu}\|_{V(D)} $.} We particularly focus on elliptic equations where the input coefficient depends affinely and non-affinely on the parameters (see Section \ref{problem_setting}). However, since our error analysis only depends on the coefficient upper bounds, we expect that the methods and results presented herein can be applied to other, more general model problems with finite parametric dimension, including nonlinear elliptic PDEs, initial value problems and parabolic equations \cite{CCS14,HS13,HS13b,HS13c}. Our approach seeks a direct estimate of { $\sum_{{\bm \nu}\in (\Lambda^{\textbf{qsi}}_M)^c}B({\bm \nu}) $} without using the Stechkin inequality. It involves a partition of { $B( (\Lambda^{\textbf{qsi}}_M)^c)$} into a family of small positive real intervals $(\mathcal{I}_j)_{j\in \mathcal{J}}$ and the corresponding splitting of { $ (\Lambda^{\textbf{qsi}}_M)^c$} into disjoint subsets $\mathcal{Q}_j$ of indices $\bm \nu$, such that $B({\bm \nu}) \in \mathcal{I}_j$. Under this process, the truncation error can be bounded as 
\begin{align*}
\sum_{{\bm \nu}\in  { (\Lambda^{\textbf{qsi}}_M)^c}}B({\bm \nu}) = \sum_{j\in \mathcal{J}} \sum_{{\bm \nu}\in \mathcal{Q}_j} B({\bm \nu}) \le \sum_{j\in \mathcal{J}} \#(\mathcal{Q}_j)\cdot \max(\mathcal{I}_j), 
\end{align*}
and therefore, the quality of the error estimate mainly depends on the approximation of cardinality of $\mathcal{Q}_j$. To tackle this problem, we develop a strategy which extends $\mathcal{Q}_j$ into continuous domain and, through relating the number of $N$-dimensional lattice points to continuous volume (Lebesgue measure), establishes a sharp estimate of the cardinality $\#(\mathcal{Q}_j)$ up to any prescribed accuracy. This development includes the utilization and extension of several results on lattice point enumeration; for a survey we refer to \cite{BR07,Gru07}. Under some weak assumptions on $B({\bm \nu})$ (which are satisfied by all existing coefficient estimates we are aware of), we achieve an asymptotic sub-exponential convergence rate of truncation error of the form $M\exp(-(\kappa M)^{1/N})$, where $\kappa$ is a constant depending on the shape and size of quasi-optimal index sets. Through several examples, we explicitly derive $\kappa$ and demonstrate the optimality of our estimate both theoretically (by proving a lower bound) and computationally (via comparison with exact calculation of truncation error). The advantage of our analysis framework is therefore twofold. First, it applies to a general class of {  coefficient decay (and correspondingly, quasi-optimal and best $M$-term approximations)}; and second, it yields sharp estimates of the asymptotic convergence rates. { For convenience, for the rest this paper, we will drop the superscript $\textbf{qsi}$ and simply refer to the quasi-optimal set of cardinality $M$ as $\Lambda_M$. }

Our paper is organized as follows. In Section \ref{problem_setting}, we describe the elliptic equations with parameterized input coefficient and necessary mathematical notations. In Section \ref{sec:holomorphy}, we present the analyticity of the solution $u$ with respect to parameter and derive coefficient estimates of Taylor and Legendre expansions of $u$. The advantage of Legendre over Taylor expansions will also be discussed. {  Our main results on the convergence analysis for a general class of multi-indexed series $\sum_{{\bm \nu}\in \mathcal{S}}B({\bm \nu})$ are established in Section \ref{anal_general_series}}. {  By means of these results, we give accessible asymptotic error estimate of several quasi-optimal and best-$M$ term polynomial approximations in Section \ref{section:conv_anal}. Finally, Section \ref{sec:validity} is devoted to further discussions on the error lower bound, as well as the pre-asymptotic estimate in a simplified case.} 
\section{Problem setting}
\label{problem_setting}
We consider solving simultaneously the following parameterized linear, elliptic PDE:
\begin{equation}
\label{elliptic_eq1}
\left\{
\begin{array}{rll}
  -\nabla\cdot (a(x,{\bm y})\nabla u(x,{\bm y})) &= f(x), &\forall (x,{\bm y})\in D\times\Gamma, \\
    u(x,{\bm y}) &= 0, &\forall (x,{\bm y})\in\partial D\times \Gamma,
\end{array}
\right.
\end{equation}  
on a bounded Lipschitz domain $D\subset \mathbb{R}^d$, with the coefficient $a(\cdot,{\bm y})$ defined on 
$\Gamma =\prod_{i=1}^N\Gamma_i \subset \mathbb{R}^N$, with $\Gamma_i =  [-1,1],\, \forall i\in \{1,\ldots,N\}$. We require the following assumption:
\begin{assumption}[Continuity and coercivity]
\label{uniform_ellipticity}
There exist constants \\
$0<a_{\min}\le a_{\max}$ such that for all $x\in \overline{D}$ and ${\bm y}\in \Gamma$
\begin{align*}
a_{\min} \le a(x,{\bm y}) \le a_{\max}.
\end{align*} 
\end{assumption} 
The Lax-Milgram lemma ensures the existence and uniqueness of solution $u$ in $V(D)\otimes L^2_{\varrho}(\Gamma)$, where $V(D) = H_0^1(D)$ and $L^2_{\varrho}(\Gamma) $ is the space of square integrable functions on $\Gamma$ with respect to the measure $\varrho({\bm y})d{\bm y}$ with $\varrho({\bm y}) = \prod_{i=1}^N \varrho_i(y_i),\, \varrho_i=\frac{1}{2},\, \forall {\bm y}\in \Gamma$. This setting represents parametric elliptic models as well as stochastic models with bounded support random coefficient. We denote $V^*(D) = H^{-1}(D)$ and, without loss of generality, assume $a_{\min} = 1$ in this work.

The corresponding weak formulation for \eqref{elliptic_eq1} is written as follows: find $u(x,{\bm y})\in V(D)\otimes L^2_{\varrho}(\Gamma)$ such that 
\begin{gather}
\label{elliptic_weak}
\begin{aligned}
&\int_\Gamma  \int_D a(x,{\bm y})\nabla u(x,{\bm y})\cdot \nabla v(x,{\bm y}) dxd{\bm y}
\\ 
&\qquad =  \int_\Gamma  \int_D  f(x)  v(x,{\bm y}) dxd{\bm y}\ \ \ \forall\, v\in V(D)\otimes L^2_{\varrho}(\Gamma).
\end{aligned}
\end{gather}

Following the arguments in \cite{CDS11}, we derive the convergence of Taylor and Legendre approximations based on the analyticity of the solution on complex domains. 
Here, the convergence is proved under the affine parameter dependence of diffusion coefficients for the Taylor series, but we relax this assumption for the Legendre series. More specifically, we only assume a holomorphic extension $a(x,{\bm z})$ of $a(x,{\bm y})$ for the complex variable ${\bm z} = (z_1,\cdots,z_N)^\top$:

\begin{assumption}[Holomorphic parameter dependence] 
\label{holomorphic_extend}
The complex continuation of $a$, represented as the map $a:\mathbb{C}^N\to L^\infty ({D})$, is a $L^\infty({D})$-valued holomorphic function on $\mathbb{C}^N$. 
\end{assumption}
This condition is easily fulfilled with $a(x,{\bm y})$ consisting of polynomials, exponential, sine and cosine functions of the variables $y_1,\cdots, y_N$. Below, we give some examples of diffusion coefficients which can be accommodated in our framework. The rigorous proofs and discussion on the advantage of Legendre over Taylor approximations will be postponed to the next section. 

\textit{Example 1.} For the input coefficient depending affinely on the parameters, i.e., 
\begin{align*}
a(x,{\bm y}) = a_0(x) + \sum_{i=1}^N y_i\psi_i(x),\ \ x\in \overline{D},\, {\bm y} \in \Gamma,
\end{align*}
where $a_0\in L^{\infty}(D)$, $(\psi_i)_{1\le i \le N} \subset L^{\infty}(D)$ such that $a$ satisfies Assumption \ref{uniform_ellipticity}; both Taylor and Legendre series approximations of $u({\bm y})$ to \eqref{elliptic_eq1} converge.  


\textit{Example 2.} Consider the input coefficient defined as
\begin{align*}
a(x,{\bm y}) = a_0(x) + \left(\sum_{i=1}^N y_i\psi_i(x)\right)^2,\ \ x\in \overline{D},\, {\bm y} \in \Gamma,
\end{align*}
with $a_0\in L^{\infty}(D),\, a_0(x)\ge a_{\min}>0\ \forall x\in \overline{D}$ and $(\psi_i)_{1\le i \le N} \subset L^{\infty}(D)$. It is easy to see that $a(x,{\bm y})$ satisfies Assumptions \ref{uniform_ellipticity}--\ref{holomorphic_extend}. Thus, the Legendre series approximation of$u({\bm y})$ to \eqref{elliptic_eq1} converges for this model.  

\textit{Example 3.} Consider the input coefficient defined as
\begin{align*}
a(x,{\bm y}) = a_0(x) + \exp \left(\sum_{i=1}^N y_i\psi_i(x)\right),\ \ x\in \overline{D},\, {\bm y} \in \Gamma,
\end{align*}
with $a_0\in L^{\infty}(D),\, a_0(x)\ge 0\ \forall x\in \overline{D}$ and $(\psi_i)_{1\le i \le N} \subset L^{\infty}(D)$.
We have $a(x,{\bm y})$ satisfies Assumptions \ref{uniform_ellipticity}--\ref{holomorphic_extend} and Legendre series approximation of $u({\bm y})$ to \eqref{elliptic_eq1} converges. 

Another framework for establishing the convergence of Legendre series was presented in \cite{CCS14} and applied to a {  more general setting} of parametric PDEs (non-elliptic, infinite dimensional noise and non-affine dependence on parameters). This approach imposes analyticity assumptions on the solution, which requires nontrivial validation in practice. Instead, in this work, we focus on elliptic equations which allows us to derive concise, minimal assumptions on the {input coefficient} (as seen above), under which the convergence of Legendre approximations holds straightforwardly. It is also worth recalling that {  our} error estimates only depend on the sharp upper bound of the polynomial coefficients. Therefore, while not studied herein, PDE models covered by \cite{CCS14,HS13,HS13b,HS13c}, bringing about same types of coefficient bounds as those considered in Section \ref{section:conv_anal}, can be treated by {  the forthcoming} analysis. 

\section{Analyticity of the solutions and estimates of the polynomial coefficients}
\label{sec:holomorphy}
Loosely speaking, the coefficients of Taylor and Legendre expansions can be estimated via three steps: 
\begin{enumerate}
\item Extending the uniform ellipticity of $a$ from $\Gamma$ to certain polydiscs/polyellipses in $\mathbb{C}^N$;
\item Proving the analyticity of the solution on those extended domains; and
\item Estimating the expansion coefficients using the analyticity properties and Cauchy's integral formula.
\end{enumerate} 

We will discuss each step in the next subsections. By $\Re(z)$ and $\Im(z)$, we denote the real and imaginary part of a complex number $z$. 

\subsection{Complex uniform ellipticity} 
The convergence of Taylor approximations is proved using the uniform ellipticity of the input coefficient in polydiscs containing $\Gamma$, based on complex analysis argument.
\begin{definition}
\label{polydisc_holomorphic}
For $0<\delta < a_{\min}$ and ${\bm \rho}$ denoting the vector $(\rho_i)_{1\le i\le N}$ with $\rho_i > 1\ \, \forall i$, we say $a(x,{\bm y})$ satisfies $(\delta,{\bm \rho})$-polydisc uniform ellipticity assumption (referred to as {\normalfont \textbf{DUE}$(\delta,{\bm \rho})$}) if there holds 
\begin{align*}
\Re(a(x,{\bm z})) \ge \delta
\end{align*} 
for all $x\in \overline{D}$ and all ${\bm z}= (z_i)_{1\le i \le N }$ contained in the polydisc 
\begin{align*}
\mathcal{O}_{\bm \rho} = \bigotimes_{1\le i\le N}\left\{ z_i\in \mathbb{C}:|z_i|\le {\rho_i}\right\}.
\end{align*}
\end{definition}

At the same time, Legendre expansions require the uniform ellipticity in smaller complex domains: the polyellipses.
\begin{definition}
\label{polyellipse_holomorphic}
For $0<\delta < a_{\min}$ and ${\bm \rho}$ denoting the vector $(\rho_i)_{1\le i\le N}$ with $\rho_i > 1\ \, \forall i$, we say $a(x,{\bm y})$ satisfies $(\delta,{\bm \rho})$-polyellipse uniform ellipticity assumption (referred to as {\normalfont \textbf{EUE}$(\delta,{\bm \rho})$}) if there holds 
\begin{align*}
\Re(a(x,{\bm z})) \ge \delta
\end{align*} 
for all $x\in \overline{D}$ and all ${\bm z}= (z_i)_{1\le i \le N }$ contained in the polyellipse 
\begin{align*}
\mathcal{E}_{\bm \rho} =
\bigotimes_{1\le i\le N}\bigg\{ z_i\in \mathbb{C}: &\Re(z_i)=\frac{\rho_i + \rho_i^{-1}}{2}\cos \phi, \\
&\left. \Im(z_i)=\frac{\rho_i - \rho_i^{-1}}{2}\sin \phi,\, \phi \in [0,2\pi) \right\}.
\end{align*}
\end{definition}

A close look at \textbf{DUE} and \textbf{EUE} reveals the advantage of Legendre over Taylor expansions. The polyellipses $\mathcal{E}_{\bm \rho}$ extend the real domain $\Gamma$ in a continuous manner, so that if ${\bm \rho}$ tends toward ${\bm 1}$, $\mathcal{E}_{\bm \rho}$ shrinks to $\Gamma$. Thus, it is hopeful that the uniform ellipticity property of $a(x,{\bm y})$ in $\Gamma$ (Assumption \ref{uniform_ellipticity}) can carry over to some polyellipses $\mathcal{E}_{\bm \rho}$ (at least with ${\bm \rho}$ close to ${\bm 1}$). In fact, we prove that \textbf{EUE} property is a consequence of Assumptions \ref{uniform_ellipticity} and \ref{holomorphic_extend}.
\begin{lemma}
\label{lemma:UE}
Let $a: \Gamma \to L^{\infty}(D)$ be a continuous function satisfying Assumptions \ref{uniform_ellipticity} and \ref{holomorphic_extend}. Then, for all $\delta< a_{\min}$, there exists a vector ${\bm\rho} = (\rho_i)_{1\le i\le N}$ with $\rho_i>1\ \forall i$ such that {\normalfont \textbf{EUE}$(\delta,{\bm \rho})$} holds.
\end{lemma}

\begin{figure}[h]
\centering
\includegraphics[height=3.3cm]{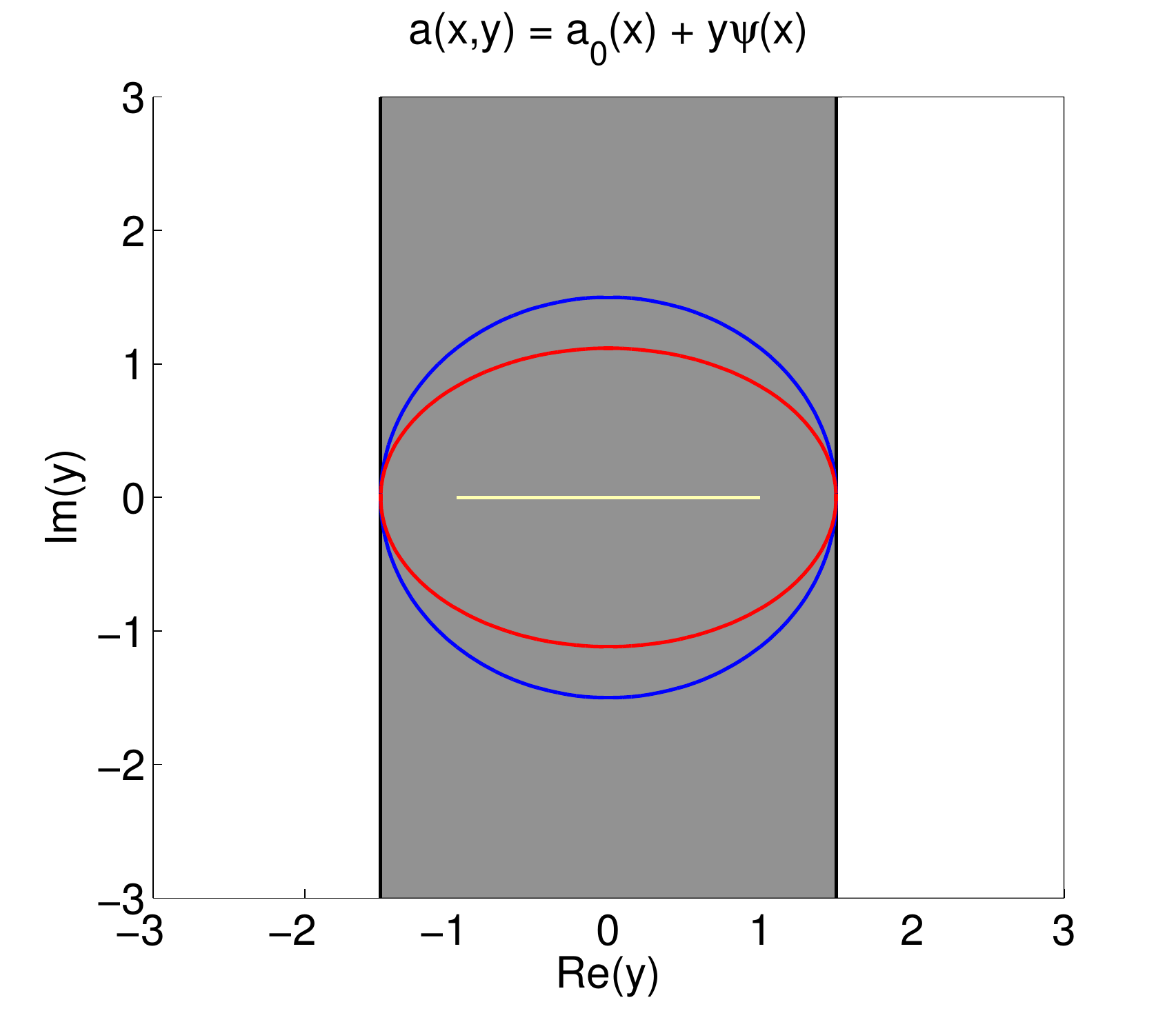} 
\includegraphics[height=3.3cm]{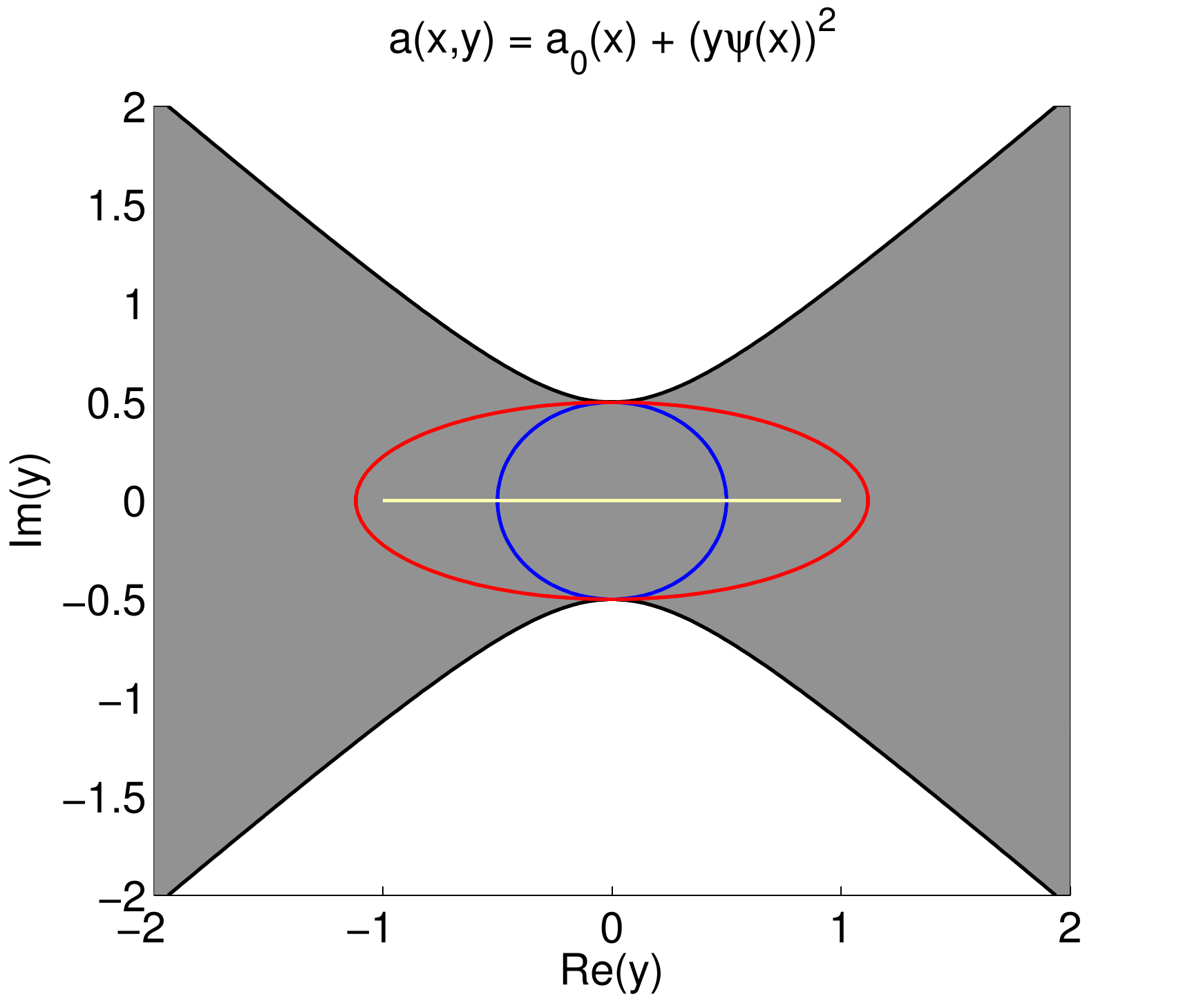} 
\includegraphics[height=3.3cm]{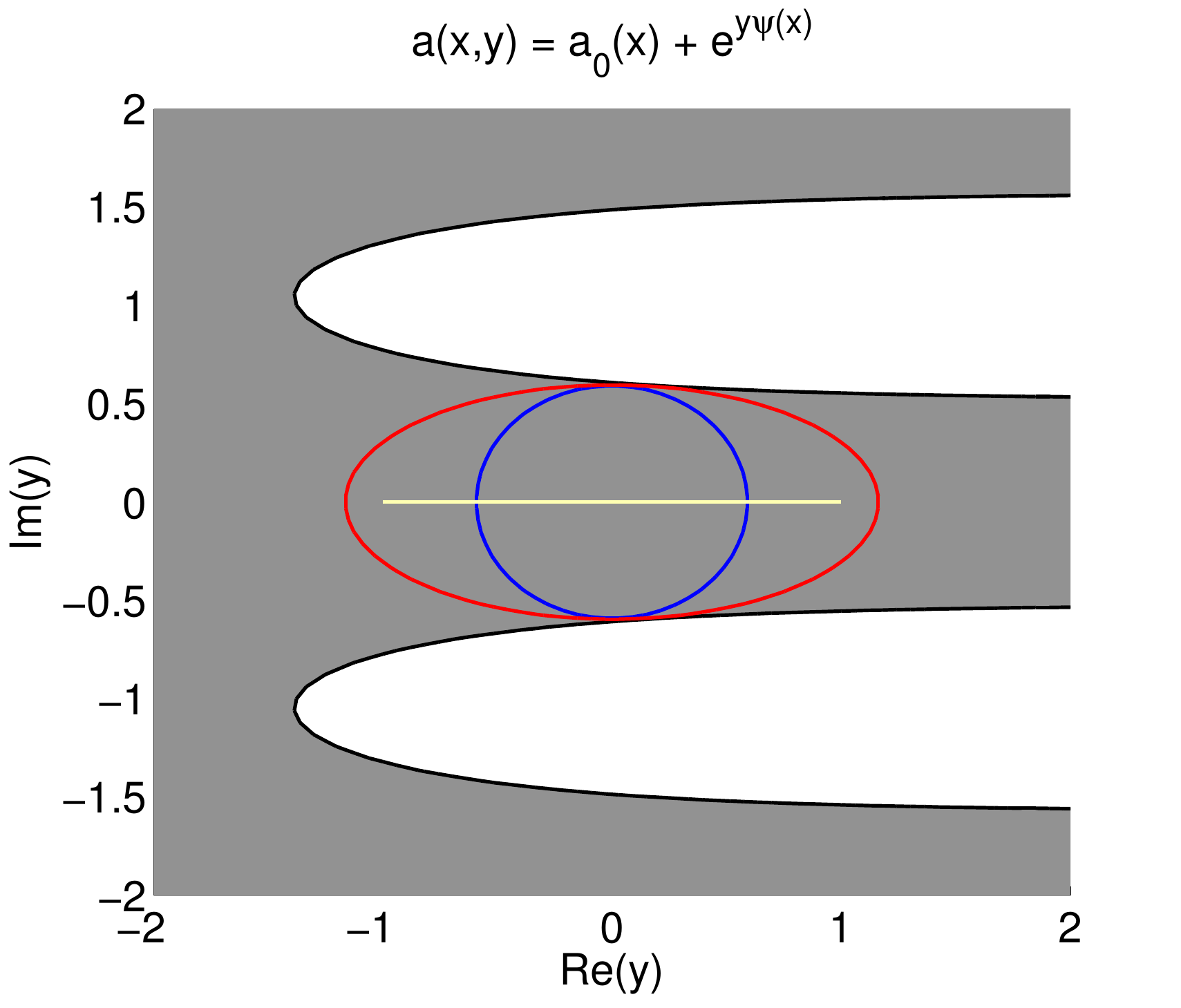} 
\caption{Domains of uniform ellipticity for some $1$-$d$ parametric input coefficients (indicated in gray). The yellow lines represent the interval $[-1,1]$. The blue and red curves are the maximal discs and ellipses which can be contained in those domains, respectively.}
\label{domUE}
\end{figure}

On the other hand, \textbf{DUE} always requires an extension of the coercive property in $\Gamma$ to the unit polydisc $\mathcal{O}_{\bm 1}$, to say the least, which is not possible generally. For illustration, the sets of $\bm z$ such that $\Re(a(x,{\bm z})) \ge \delta$ for all $x\in \overline{D}$ with some fixed $\delta>0$ (referred to as the domains of uniform ellipticity) are plotted in Figure \ref{domUE} for some typical $1$-dimensional parametric coefficients. The maximal ellipses and discs contained in these domains are shown. We observe that for the affine coefficient, the set spans unrestrictedly along the imaginery axis, and discs covering $\Gamma$ can easily be placed inside. It highlights the success of Taylor approximations for parameterized models which depend affinely on the parameters, \cite{CDS11,HS13,HS13b}. This property however no longer holds for non-affine, yet holomorphic diffusion coefficients. Taylor approximations for these cases can be treated by a real analysis approach, but under additional strong constraints, see \cite{BTNT12}. 

We close this subsection with a proof of Lemma \ref{lemma:UE}. 

\sproof
(of Lemma \ref{lemma:UE}).
Since $a({\bm z})$ is holomorphic in $\mathbb{C}^N$, we have $\Re(a({\bm z}))$ is a continuous mapping. By Heine-Cantor theorem, $\Re(a({\bm z}))$ is uniformly continuous in any compact subset of $\mathbb{C}^N$. Fixing a $0<\delta < a_{\min}$, without loss of generality, we can choose $\xi >0$ such that $\forall {\bm z} \in \mathbb{C}^N,\, \forall {\bm z}' \in \Gamma$ satisfying $\|{\bm z} - {\bm z}'\|\le \xi$, there holds
\begin{align*}
\left \|\Re(a({\bm z})) - \Re(a({\bm z}')) \right \|_{L^{\infty}(D)} \le a_{\min} -\delta.
\end{align*}
This implies 
\begin{align*}
\Re(a(x,{\bm z})) \ge \delta - a_{\min} + \Re(a(x,{\bm z}')) \ge \delta,
\end{align*}
for all $x\in \overline{D},\, {\bm z}\in \mathbb{C}^N$ such that $\|{\bm z} - {\bm z}'\|\le \xi$ with some $ {\bm z}' \in \Gamma$. Denoting $ \Gamma_\xi = \{ {\bm z}\in \mathbb{C}^N: \mbox{dist}({\bm z},\Gamma) \le \xi \} $, we proceed to prove there exists ${\bm \rho} = (\rho_i)_{1\le i\le N}$ with $\rho_i>1\ \forall i$ such that the polyellipse $\mathcal{E}_{\bm \rho}$ is included in $\Gamma_\xi$. 

First, consider the ``polyrectangle'' 
\begin{align*}
\Xi =  \bigotimes_{1\le i\le N}\left\{ z_i\in \mathbb{C}:|\Re(z_i)|\le 1 + \frac{\xi}{\sqrt{2N}},\, |\Im(z_i)| \le \frac{\xi}{\sqrt{2N}}\right\}, 
\end{align*}
we will show that $\Xi\subset \Gamma_{\xi}$. Indeed, for every ${\bm z}\in \Xi$, choose ${\bm z}'= (z'_i)_{1\le i\le N}$ as follows: if $|\Re(z_i)|\le 1$, then $z'_i = \Re(z_i)$; otherwise, $z'_i = \mbox{sgn}(\Re(z_i))$. It is easy to see that ${\bm z}'\in \Gamma$. Furthermore, for all $ i \in \{1,\ldots, N\}$,
\begin{align*}
|\Re(z_i) - \Re(z'_i)| \le \frac{\xi}{\sqrt{2N}}, \ \ |\Im(z_i) - \Im(z'_i)| \le \frac{\xi}{\sqrt{2N}}.
\end{align*}
Thus, $|z_i - z'_i| \le \frac{\xi}{\sqrt{N}} $ and we have
\begin{align*}
\|{\bm z} - {\bm z}' \| \le \left(\sum_{i=1}^N |z_i - z'_i|^2 \right)^{1/2} \le \xi, 
\end{align*}
This gives ${\bm z}\in \Gamma_{\xi}$ and $\Xi\subset \Gamma_{\xi}$. 

It remains to find ${\bm \rho}$ satisfying $\mathcal{E}_{\bm \rho} \subset \Xi$. To make this hold, we only need to select ${\bm \rho}$ such that the lengths of axes of each ellipse are less than the lengths of corresponding sizes of the rectangle, i.e.,   
\begin{align*}
\frac{\rho_i + \rho_i^{-1}}{2} \le 1 + \frac{\xi}{\sqrt{2N}},\ \ \mbox{ and }\ \ \frac{\rho_i - \rho_i^{-1}}{2} \le  \frac{\xi}{\sqrt{2N}}.
\end{align*}

The choice of $\rho_i = \frac{\xi}{\sqrt{2N}} + \sqrt{\frac{\xi^2}{2N} + 1} > 1$ fulfills this condition, with which $\mathcal{E}_{\bm \rho} \subset \Xi \subset \Gamma_{\xi}$. There follows 
\begin{align*}
\Re(a(x,{\bm z})) \ge \delta
\end{align*} 
for all $x\in \overline{D},\, {\bm z}\in \mathcal{E}_{\bm \rho}$ and $a$ satisfies  {\normalfont \textbf{EUE}$(\delta,{\bm \rho})$}, as desired. 
\fproof

\subsection{Analyticity of the solutions with respect to the parameters} $ $
If \textbf{DUE}/\textbf{EUE} holds, according to the Lax-Milgram theorem, $u({\bm z})\in V(D)$ is defined and uniformly bounded in certain polydiscs $\mathcal{O}_{\bm \rho}$/polyellipses $\mathcal{E}_{\bm \rho}$ containing $\Gamma$. Exploiting this fact and the analyticity of $a({\bm z})$ in $\mathbb{C}^N$, we establish the analyticity of the map ${\bm z} \mapsto u({\bm z})$. The results given in this section { essentially follow those in Section 2.1 of \cite{CDS11}, but apply to the more general cases} of smooth, non-affine diffusion coefficients. 




\begin{theorem}
\label{holomorphic}
Assume that the coefficient $a(x,{\bm y})$ satisfies Assumptions \ref{uniform_ellipticity}--\ref{holomorphic_extend}. If {\normalfont \textbf{DUE}$(\delta,{\bm \rho})$} ({\normalfont \textbf{EUE}$(\delta,{\bm \rho})$} correspondingly) holds for some $0<\delta <a_{\min}$ and ${\bm \rho}=(\rho_i)_{1\le i\le N}$ with $\rho_i > 1\ \, \forall i$, then the function ${\bm z}\mapsto u({\bm z})$ is holomorphic in an open neighborhood of the polydisc $\mathcal{O}_{{\bm \rho}}$ (the polyellipse $\mathcal{E}_{{\bm \rho}}$ correspondingly). 
\end{theorem}

\sproof
We will prove this theorem for $a(x,{\bm y})$ satisfying \textbf{DUE}$(\delta,{\bm \rho})$. The other case should follow similarly. Defining 
\begin{align*}
\mathcal{A} = \left\{{\bm z}\in \mathbb{C}^N: \Re(a(x,{\bm z})) > \frac{\delta}{2}\mbox{ for all }x\in \overline{D} \right\},
\end{align*}
the proof consists of two steps showing that 
\begin{enumerate}
\item $\mbox{int}(\mathcal{A})$ is an open neighborhood of the polydisc $\mathcal{O}_{\bm \rho}$; and 
\item the map ${\bm z}\mapsto u({\bm z})$ is holomorphic in $\mbox{int}(\mathcal{A})$.
\end{enumerate}
Here, $\mbox{int}(\mathcal{A})$ is the interior of $\mathcal{A}$. 

First, let us choose an arbitrary element $\tilde{\bm z}$ in $\mathcal{O}_{\bm \rho}$. For $B({\bm z},r)$, we denote the open ball radius $r$ centered at ${\bm z}$ in $\mathbb{C}^N$. Observing that the map ${\bm z}\mapsto a({\bm z})$ is holomorphic in $\mathbb{C}^N$, we have ${\bm z}\mapsto \Re(a({\bm z}))$ is a continuous function in $\mathbb{C}^N$. There exists $r_{\tilde{\bm z}}>0$ depending on $\tilde{\bm z}$ such that for all ${\bm z}\in B(\tilde{\bm z},r_{\tilde{\bm z}}) $,
\begin{align*}
\left \|\Re(a(\tilde {\bm z})) - \Re(a({\bm z}))\right \|_{L^{\infty}(D)} < \frac{\delta}{2}. 
\end{align*}
This gives 
\begin{align*}
\Re(a(x,  {\bm z})) > \Re(a(x, \tilde {\bm z})) -\frac{\delta}{2} \ge \frac{\delta}{2},\ \ \ \forall x\in \overline{D},\, {\bm z} \in B(\tilde{\bm z},r_{\tilde{\bm z}}),
\end{align*}
and $B(\tilde{\bm z},r_{\tilde{\bm z}})\subset \mathcal{A}$ for all $ \tilde{\bm z}\in \mathcal{O}_{\bm \rho}$. We obtain $\tilde{\bm z}\in \mbox{int}(\mathcal{A})$ for all $\tilde{\bm z}\in \mathcal{O}_{\bm \rho}$, which concludes Step 1. 

{  As $a({\bm z})$ is holomorphic in $\mathbb{C}^N$, for all ${\bm z}\in \mbox{int}(\mathcal{A})$, there exists
\begin{align*}
\partial_i a({\bm z}) := \lim_{h\to 0} \frac{a({\bm z} + h{\bm e}_i)-a({\bm z})}{h} \in L^{\infty} (D), 
\end{align*}
where $i\in\{1,\ldots,N\}$ and ${\bm e}_i$ denotes the Kronecker sequence with $1$ at index $i$ and $0$ at other indices. The proof of Step 2 is then similar to that of Lemma 2.2, \cite{CDS11} and would be omitted here. 

}
\fproof


\subsection{Estimates of the polynomial coefficients}
\label{estim_coef}
Under the analyticity properties established in Theorem \ref{holomorphic}, the convergence of Taylor and Legendre expansions of the solutions, as well as estimates of the expansion coefficients, are well-studied, e.g., in \cite{BNTT14,CCS14,CDS11}. In this subsection, we {  review without proof} those results in the context of finite-dimensional, possibly non-affine parametric coefficients. Recall that $\mathcal{S} = \{{\bm \nu} = (\nu_i)_{1\le i\le N} : \nu_i\in \mathbb{N}\}$. For all ${\bm \nu}\in\mathcal{S}$, we introduce the multivariate notations $|{\bm \nu}| := \sum_{1\le i\le N} \nu_i$, ${\bm \nu} ! := \prod_{1\le i\le N} \nu_i !$ and define the partial derivative $\partial^{\bm \nu} u := \frac{\partial^{|{\bm \nu}| }u}{\partial^{\nu_1} z_1\ldots \partial^{\nu_N} z_N}$. The Taylor series of $u({\bm y})$ reads 
\begin{align}
\label{def:taylor}
u({\bm y}) = \sum_{{\bm \nu}\in \mathcal{S}} t_{{\bm \nu}}{\bm y}^{\bm \nu},
\end{align}
where the coefficients $t_{\bm \nu} \in V(D)$ are defined as 
\begin{align*}
t_{\bm \nu} := \frac{1}{{\bm \nu}!}\partial^{\bm \nu} u(0),\ \  {\bm \nu}\in \mathcal{S}. 
\end{align*}
The convergence of the Taylor expansion in $\Gamma$ and estimates of $\|t_{\bm \nu}\|_{V(D)}$ are given in the following.
\begin{proposition}
\label{theorem:taylor_est}
Assume that the coefficient $a(x,{\bm y})$ satisfies Assumptions \ref{uniform_ellipticity}--\ref{holomorphic_extend}. If {\normalfont \textbf{DUE}$(\delta,{\bm \rho})$} holds for some $0<\delta <a_{\min}$ and ${\bm \rho}=(\rho_i)_{1\le i\le N}$ with $\rho_i > 1\ \, \forall i$ then the Taylor series $\sum_{{\bm \nu}\in \mathcal{S}} t_{{\bm \nu}}{\bm y}^{\bm \nu}$ converges uniformly towards $u({\bm y})$ in $\Gamma$. Furthermore, we have the estimate 
\begin{align}
\|t_{\bm \nu}\|_{V(D)} \le \frac{\|f\|_{V^*(D)}}{\delta} {\bm \rho}^{-{\bm \nu}}.  \label{taylor_est}
\end{align}
\end{proposition}

{  [Proof was removed.]}

On the other hand, the tensorized Legendre series of $u({\bm y})$ is defined as
\begin{align}
\label{def:legendre1}
u({\bm y}) = \sum_{{\bm \nu}\in \mathcal{S}} u_{{\bm \nu}}P_{\bm \nu}({\bm y}),
\end{align}
where $P_{\bm \nu}({\bm y}) = \prod_{i=1}^N P_{\nu_i}(y_i)$, with $P_{\nu_i}$ denoting the monodimensional Legendre polynomials of degree $\nu_i$ according to $L^{\infty}$ normalization $\|P_{\nu_i}\|_{L^{\infty}([-1,1])} = P_{\nu_i}(1) = 1$.

A second type of Legendre expansion, which employs the $L^2$ normalized version of $P_{\bm \nu}$, is also considered. Denote the multivariate polynomials $ L_{\bm \nu}({\bm y}) = \prod_{i=1}^N L_{\nu_i}(y_i)$, with $L_{\nu_i}(y_i)$ given by
\begin{align*}
L_{\nu_i}(y_i) := \sqrt{2\nu_i + 1} P_{\nu_i}(y_i).
\end{align*}
The Legendre series in this case can be written as 
\begin{align}
\label{def:legendre2}
u({\bm y}) = \sum_{{\bm \nu}\in \mathcal{S}} v_{{\bm \nu}}L_{\bm \nu}({\bm y}).
\end{align}
We note that the coefficients $u_{\bm \nu},\, v_{\bm \nu}\in V(D)$ are defined by
\begin{align}
\label{leg_coef_relation}
v_{\bm \nu} = \int_\Gamma u({\bm y}) L_{\bm \nu}({\bm y}) \varrho({\bm y})d{\bm y}\ \ \mbox{ and }\ \ 
u_{\bm \nu} = {v_{\bm \nu}}{\left ( \prod_{i=1}^N{(2\nu_i + 1)} \right )^{1/2}}.
\end{align} 
The following proposition establishes estimates of $\|u_{\bm \nu}\|_{V(D)},\, \|v_{\bm \nu}\|_{V(D)}$ and the convergence of the Legendre expansions of $u$ in $\Gamma$. 

\begin{proposition}
\label{theorem:legendre_est}
Assume that the coefficient $a(x,{\bm y})$ satisfies Assumptions \ref{uniform_ellipticity}--\ref{holomorphic_extend}. If {\normalfont \textbf{EUE}$(\delta,{\bm \rho})$} holds for some $0<\delta <a_{\min}$ and ${\bm \rho}=(\rho_i)_{1\le i\le N}$ with $\rho_i > 1\ \, \forall i$ then we have the estimates 
\begin{gather}
 \label{legendre_est}
\begin{aligned}
\|u_{\bm \nu}\|_{V(D)} \le C_{{\bm \rho},\delta} {\bm \rho}^{-{\bm \nu}}\prod_{i=1}^N {(2\nu_i+1)},\ \ \|v_{\bm \nu}\|_{V(D)} \le C_{{\bm \rho},\delta} {\bm \rho}^{-{\bm \nu}}\prod_{i=1}^N\sqrt{2\nu_i+1},
\end{aligned}
\end{gather}
where $C_{{\bm \rho},\delta} = \frac{\|f\|_{V^*(D)}}{\delta} \prod_{i=1}^{N}\frac{\ell(\mathcal{E}_{\rho_i})}{4(\rho_i -1)}$ with $\ell(\mathcal{E}_{\rho_i})$ denoting the perimeter of the ellipse $\mathcal{E}_{\rho_i}$.

Consequently, the Legendre series $\sum_{{\bm \nu}\in \mathcal{S}} u_{{\bm \nu}}P_{\bm \nu}$ and $\sum_{{\bm \nu}\in \mathcal{S}} v_{{\bm \nu}}L_{\bm \nu}$ converge towards $u$ in $L^{\infty}(\Gamma,V(D))$. The series $\sum_{{\bm \nu}\in \mathcal{S}} v_{{\bm \nu}}L_{\bm \nu}$ also converges towards $u$ in $V(D)\otimes L^2_{\varrho}(\Gamma)$.
\end{proposition}

{  [Proof was removed.]}

Under Assumptions \ref{uniform_ellipticity}--\ref{holomorphic_extend}, we remark that \textbf{EUE} and, if adding affine dependence on parameters, \textbf{DUE} normally hold for infinitely many couples of $(\delta,{\bm \rho})$. We call the set of all $(\delta,{\bm \rho})$ such that \textbf{EUE}$(\delta,{\bm \rho})$/\textbf{DUE}$(\delta,{\bm \rho})$ is fulfilled the admissible set and denote it by $\bm{Ad}$ for both cases. For a fixed ${\bm \nu} \in \mathcal{S}$, the best coefficient bounds given by Propositions \ref{theorem:taylor_est} and \ref{theorem:legendre_est} will be 
\begin{align}
\label{bound_opt}
\|t_{\bm \nu}\|_{V(D)} \le \inf\limits_{(\delta,{\bm \rho})\in \bm{Ad}} \frac{\|f\|_{V^*(D)}}{\delta} {\bm \rho}^{-{\bm \nu}},\ \  \|u_{\bm \nu}\|_{V(D)} \le  \inf\limits_{(\delta,{\bm \rho})\in \bm{Ad}}  C_{{\bm \rho},\delta} {\bm \rho}^{-{\bm \nu}}\prod_{i=1}^N {(2\nu_i+1)}. 
\end{align}
Finding an efficient computation of these infimums and algorithm to construct the corresponding quasi-optimal index sets is an open question. In the specific case where the basis functions $\psi_i$ have non-overlapping supports, however, the vectors $\bm \rho$ solving the minimization problems in \eqref{bound_opt} can be found easily. In this case, the best a priori estimates retrieve the forms \eqref{taylor_est} and \eqref{legendre_est}. Recent studies have shown that although these theoretical bounds are not sharp, they construct {  quite accurate} polynomial spaces, see \cite{BNTT14}. 

\section{Asymptotic convergence analysis for a general class of multi-indexed series}
\label{anal_general_series}

{  Consider a multi-indexed sequence of coefficient estimates written in the form $(B({\bm \nu}))_{{\bm \nu}\in \mathcal{S}} \equiv (e^{-b({\bm \nu})})_{{\bm \nu}\in \mathcal{S}}$. In this section, we introduce a new, generalized approach to estimating the asymptotic convergence of $\sum_{{\bm \nu}\in \Lambda_M^c} e^{-b({\bm \nu})}$ with respect to $M$, under some general assumptions of $b$ which accommodate most types of Taylor and Legendre coefficient bounds established in current literature. We recall that this truncation error represents the error of quasi-optimal methods, as well as accessible error of best $M$-term approximations. It is enough to conduct the analysis with $\Lambda_M^c$ being the sets of all ${\bm \nu}$ such that $e^{-b({\bm \nu})} < e^{-J}$ with some $J\in \mathbb{N}$. }

Our method can be summarized as follows. First, we split $\Lambda^c_M$ into a family $(\mathcal{Q}_j)_{j\in \mathbb{N},\, j\ge J}$ of disjoint subsets of ${\mathcal{S}}$ based on values of $e^{-b({\bm \nu})}$, where $\mathcal{Q}_j$ contains ${\bm \nu}$ satisfying $e^{-j-1} \le e^{-b({\bm \nu})} < e^{-j}$, so that the truncation error can be bounded as 
\begin{align}
\label{sec4:intro}
\sum_{{\bm \nu}\in \Lambda^c_M}e^{-b({\bm \nu})} = \sum_{j \ge J} \sum_{{\bm \nu}\in \mathcal{Q}_j} e^{-b({\bm \nu})} \le \sum_{j\ge {J}} \#(\mathcal{Q}_j)\cdot e^{-j}.
\end{align}
Obviously, finding a sharp approximation of $\#(\mathcal{Q}_j)$ is central to estimate \eqref{sec4:intro}. We define the \textit{superlevel sets} $\mathcal{P}_j$ of $N$-dimensional real points 
\begin{align}
\mathcal{P}_j := \{{\bm \nu} \in [0,\infty)^N: e^{-b({\bm \nu})} \ge e^{-j}\} = \{{\bm \nu}\in [0,\infty)^N: b({\bm \nu})\le j\}, \label{define:suplevel_set}
\end{align} 
and, with notice that $\#(\mathcal{Q}_j) = \#(\mathcal{P}_{j+1} \cap \mathbb{Z}^N) - \#(\mathcal{P}_{j} \cap \mathbb{Z}^N)$, seek to count points with integer coordinates in $\mathcal{P}_j$. An appealing approach to solving this problem is to study the interplay between $\#(\mathcal{P}_{j}\cap \mathbb{Z}^N)$ and the continuous volume (Lebesgue measure) of $\mathcal{P}_j$. We first employ the following well-known result in measure theory, reflecting the intuitive fact that for a geometric body $\mathcal{P}$ in $\mathbb{R}^N$, the volume of $\mathcal{P}$, denoted by $|\mathcal{P}|$, can be approximated by the number of shrunken integer points inside $\mathcal{P}$, see, e.g., Section 7.2 in \cite{Gru07} and Section 1.1 in \cite{Tao11}.

\begin{lemma}
\label{point_count}
Suppose $\mathcal{P} \subset \mathbb{R}^N$ is a bounded Jordan measurable set. For $j \in \mathbb{N},\, j>0$, there holds 
\begin{align}
\label{eq:point_count}
|\mathcal{P}|=  \lim_{j \to  \infty} \frac{1}{j^N}\cdot \#(\mathcal{P}\cap \frac{1}{j}\, \mathbb{Z}^N)  = \lim_{j \to \infty} \frac{1}{j^N}\cdot \#(j\mathcal{P}\cap \mathbb{Z}^N) .
\end{align}
\end{lemma}

Concerning our goal of estimating \eqref{sec4:intro}, Lemma \ref{point_count} has an interesting consequence: If $b({\bm \nu})$ is defined such that $\frac{1}{j} \mathcal{P}_j = \mathcal{P}$, $\forall j\in \mathbb{N}$, with some $\mathcal{P}\subset\mathbb{R}^N$, one obtains a simple asymptotic formula for $\#(\mathcal{P}_j\cap \mathbb{Z}^N)$:
\begin{align}
\#(\mathcal{P}_j\cap \mathbb{Z}^N) \simeq j^N |\mathcal{P}|. \label{core_est}
\end{align} 
Such approximation is powerful since, loosely speaking, it would allow replacing $\#(\mathcal{Q}_j) $ by $((j+1)^N -  j^N) |\mathcal{P}|$ and reduce \eqref{sec4:intro} to a much easier, yet equivalent problem of estimating the truncation error via
\begin{align}
\label{sec4:intro2}
\sum_{{\bm \nu}\in \Lambda^c_M}e^{-b({\bm \nu})} \lesssim \sum_{j\ge J} |\mathcal{P}| ((j+1)^{N} - j^N)e^{-j}.
\end{align}
The property that the sets $\frac{1}{j}\mathcal{P}_j$ are unchanged over $j\in \mathbb{N}$ is, however, restrictive, corresponding to only a few types of coefficient upper bounds, for instance, $b({\bm \nu})$ is linear in ${\bm \nu}$. For this approach of estimation to be considered useful in general quasi-optimal approximation setting, this condition needs to be relaxed. 


For the technicality, we now extend definition \eqref{define:suplevel_set} to equip the superlevel sets with {real} indices: for $\tau\in (0,\infty)$, define
\begin{align}
\mathcal{P}_\tau := \{{\bm \nu} \in [0,\infty)^N: e^{-b({\bm \nu})} \ge e^{-\tau}\} = \{{\bm \nu}\in [0,\infty)^N: b({\bm \nu})\le \tau\}. \label{define:suplevel_set2}
\end{align} 
Note that the assertion of Lemma \ref{point_count} still holds if replacing $j\in \mathbb{N}$ by $\tau \in (0,\infty)$.    
We establish, in Lemma \ref{point_count2} below, formula \eqref{core_est} under some weaker assumptions on $(\mathcal{P}_\tau)_{\tau\in\mathbb{R}^+}$:  
\begin{enumerate}[\ \ \ i)] 
\item $\mathcal{P}_\tau$ is Jordan measurable for countably infinite $\tau\in (0,\infty)$, 
\item The chain $(\frac{1}{\tau}\mathcal{P}_\tau)_{\tau\in\mathbb{R}^+}$ is either ascending or descending towards a Jordan measurable \textit{limiting set} $\mathcal{P}\subset \mathbb{R}^N$ with $0<|\mathcal{P}|<\infty$. 
\end{enumerate} 
As we shall see later, these properties are satisfied by most existing polynomial coefficient estimates.


\begin{lemma}
\label{point_count2}
Suppose $(\mathcal{P}_\tau)_{\tau\in \mathbb{R}^+}$ is a family of bounded Lebesgue measurable sets in $\mathbb{R}^N$ satisfying either 
{\allowdisplaybreaks
\begin{align}
\label{subset}
\frac{1}{\tau_1}\mathcal{P}_{\tau_1} &\subset \frac{1}{\tau_2}\mathcal{P}_{\tau_2},\ \forall \tau_1\ge \tau_2 > 0,
\\
\label{supset}
\text{or}\ \ \ \ \ \ \frac{1}{\tau_1}\mathcal{P}_{\tau_1} &\supset \frac{1}{\tau_2}\mathcal{P}_{\tau_2},\ \forall \tau_1\ge \tau_2 > 0.
\end{align} 
}
Denote $\mathcal{P} = \bigcap\limits_{\tau\in\mathbb{R}^+}\frac{1}{\tau}\mathcal{P}_\tau$ if \eqref{subset} holds and $\mathcal{P} = \bigcup\limits_{\tau\in\mathbb{R}^+}\frac{1}{\tau}\mathcal{P}_\tau$ for the other case. If $\mathcal{P}$ is bounded Jordan measurable, $|\mathcal{P}|>0$, and there exists a sequence $(\tau_j)_{j\in\mathbb{N}}$ with $\tau_j\to\infty$ such that $\mathcal{P}_{\tau_j}$ is Jordan measurable for all $j$, there follows 
\begin{align}
|\mathcal{P}| = \lim_{\tau\to \infty} \frac{1}{\tau^{ N}}\cdot \#(\mathcal{P}_\tau\cap \mathbb{Z}^N). \label{eq:5b}
\end{align}
\end{lemma}

\sproof
We will give a proof with $(\mathcal{P}_\tau)_{\tau\in \mathbb{R}^+}$ satisfying \eqref{subset}. The other case can be shown analogously. 
Let $\varepsilon$ be an arbitrary positive number. By Lemma \ref{point_count},
\begin{align*} 
\frac{1}{\tau^{ N}}\cdot\# (\tau \mathcal{P} \cap \mathbb{Z}^N) \to |\mathcal{P}| \ \mbox{ as }\ \tau\to\infty.
\end{align*} 
Since $\mathcal{P}\subset \frac{1}{\tau}\mathcal{P}_\tau\  \forall \tau$, we can choose $T_1 > 0$ such that $\forall \tau>T_1$,
\begin{align}
\label{eq:6}
 |\mathcal{P}| -\varepsilon \le \frac{1}{\tau^{ N}}\cdot \#(\tau \mathcal{P}\cap \mathbb{Z}^N) \le \frac{1}{\tau^{ N}} \cdot \#(\mathcal{P}_\tau\cap \mathbb{Z}^N).
\end{align}
On the other hand, from $\mathcal{P} = \bigcap_{\tau\in\mathbb{R}^+}\frac{1}{\tau}\mathcal{P}_\tau$, it yields $|\mathcal{P}| = \lim\limits_{\tau\to \infty} \left | \frac{1}{\tau} \mathcal{P}_\tau\right | $. Let us pick an $L>0$ so that $\mathcal{P}_L$ is Jordan measurable and $\left | \frac{1}{L} \mathcal{P}_{L}\right | \le |\mathcal{P}| + \frac{\varepsilon}{2}$. By Lemma \ref{point_count},
\begin{align*}
\frac{1}{\tau^{ N}}\cdot\# \left(\frac{\tau}{L} \mathcal{P}_L \cap \mathbb{Z}^N\right) \to \left |\frac{\mathcal{P}_L}{L}\right |\ \mbox{ as } \tau\to\infty.
\end{align*}
There exists $T_2>L$ satisfying $\forall \tau>T_2$,
\begin{align*}
\frac{1}{\tau^{ N}}\cdot\# \left(\frac{\tau}{L} \mathcal{P}_L \cap \mathbb{Z}^N\right) \le \left |\frac{\mathcal{P}_L}{L} \right | + \frac{\varepsilon}{2} \le |\mathcal{P}| + \varepsilon. 
\end{align*}
Since $\tau> L$, we have $\mathcal{P}_\tau \subset \frac{\tau}{L}\mathcal{P}_L$, which gives 
\begin{align}
\label{eq:6b}
\frac{1}{\tau^{ N}}\cdot\# \left( \mathcal{P}_\tau \cap \mathbb{Z}^N\right) \le \frac{1}{\tau^{ N}}\cdot\# \left(\frac{\tau}{L} \mathcal{P}_L \cap \mathbb{Z}^N\right) \le |\mathcal{P}| + \varepsilon. 
\end{align}
Combining \eqref{eq:6} and \eqref{eq:6b} proves \eqref{eq:5b}.
\fproof

Lemma \ref{point_count2} provides us with an asymptotic formula of the form \eqref{core_est} to approximate the number of integer points inside $\mathcal{P}_\tau$, under some conditions on $(\mathcal{P}_\tau)_{\tau\in\mathbb{R}^+}$. Given a coefficient upper bound $e^{-b({\bm \nu})}$, it is desirable to derive properties of $b({\bm \nu})$ such that its corresponding superlevel sets $(\mathcal{P}_\tau)_{\tau\in\mathbb{R}^+}$ fulfill these conditions. For all ${\bm \nu} \in [0,\infty)^N$, define the map $H_{\bm \nu}: (0,\infty) \to \mathbb{R}$ as
\begin{align*}
H_{\bm \nu}(\tau) = \frac{1}{\tau}b(\tau{\bm \nu}),\, \forall \tau \in (0,\infty). 
\end{align*}   
We proceed to state and validate the following assumptions on ${b({\bm \nu})}$. 

\begin{assumption}\label{sec4:assump1}
The map $b:[0,\infty)^N \to\mathbb{R}$ satisfies 
\begin{enumerate}[\ \ \ \ \ \ 1.] 
\item $b({\bm 0}) = 0$ and $b$ is continuous in $[0,\infty)^N$,
\item $H_{\bm \nu}$ is either increasing in $(0,\infty)$ for all ${\bm \nu}\in [0,\infty)^N $ or decreasing in $(0,\infty)$ for all ${\bm \nu}\in [0,\infty)^N$,
\item $b({\bm \nu}) \in \Theta(|{\bm \nu}|)$. In other words, there exists $0< c < C$ such that $c|{\bm \nu}| < b({\bm \nu}) < C|{\bm \nu}|$ as ${\bm \nu} \to \infty$. 
\end{enumerate}
\end{assumption}

\begin{lemma}
\label{lemma:prop.of.b}
Assume that $b:[0,\infty)^N \to\mathbb{R}$ satisfies Assumption \ref{sec4:assump1}. For $\tau\in (0,\infty)$, denote $\mathcal{P}_\tau = \left\{\bm{\nu} \in [0,\infty)^N:\, b(\bm{\nu})\le \tau \right\}$. Let 
\begin{equation}
\label{define.P}
 \mathcal{P}= \left\{
\begin{array}{rll}
  \bigcap\limits_{\tau\in\mathbb{R}^+}\left(\frac{1}{\tau}\!\mathcal{P}_\tau\right)\mbox{ if }H_{\bm \nu}\mbox{ is increasing }\, \forall {\bm \nu}\in [0,\infty)^N, 
 \\
\bigcup\limits_{\tau\in\mathbb{R}^+}\left(\frac{1}{\tau}\!\mathcal{P}_\tau\right)\mbox{ if }H_{\bm \nu}\mbox{ is decreasing }\, \forall {\bm \nu}\in [0,\infty)^N.
\end{array}
\right.
\end{equation}
Then, $0<|\mathcal{P}|<\infty$. If $\mathcal{P}$ is Jordan measurable, there holds
\begin{align}
|\mathcal{P}| = \lim_{\tau\to \infty} \frac{1}{\tau^{ N}}\cdot \#(\mathcal{P}_\tau\cap \mathbb{Z}^N). \label{eq:6c}
\end{align}
\end{lemma}

\sproof
From the continuity of $b$ in $[0,\infty)^N$ (Assumption \ref{sec4:assump1}.1), $\mathcal{P}_\tau$ is Jordan measurable for all except a countable number of values of $\tau$ (see \cite{Fri33}). 

Next, from Assumption \ref{sec4:assump1}.2, if $H_{\bm \nu}$ is increasing for all $ {\bm \nu}$, one has $\frac{1}{\tau_2}b(\tau_2{\bm \nu}) \le \frac{1}{\tau_1}b(\tau_1{\bm \nu}) \le 1,\ \forall \tau_1 \ge \tau_2 > 0,\, \forall {\bm \nu}\in [0,\infty)^N$, which implies 
\begin{align*}
\frac{1}{\tau_1}\mathcal{P}_{\tau_1}\subset \frac{1}{\tau_2}\mathcal{P}_{\tau_2},\ \forall \tau_1 \ge \tau_2 > 0. 
\end{align*}

Since $b(\bm{\nu})$ converges towards $+\infty $ as $\bm{\nu} \to \infty$, $\mathcal{P}_\tau$ is bounded for every $\tau\in  (0,\infty)$. It is trivial that $\mathcal{P}=\bigcap_{\tau\in \mathbb{R}^+}\left(\frac{1}{\tau}\!\mathcal{P}_\tau\right)$ is bounded. Let ${\bm \nu}\notin \mathcal{P}$, we have ${\bm \nu} \notin \frac{1}{\tau}\!\mathcal{P}_\tau $ for $\tau$ large enough. Combining with Assumption \ref{sec4:assump1}.3 yields $C \tau |{\bm \nu}| > b(\tau {\bm \nu}) > \tau$. Thus, $B({\bm 0},1/C)\subset \mathcal{P}$ and $|\mathcal{P}| > 0$. 

If, on the other hand, $H_{\bm \nu}$ is decreasing for all ${\bm \nu}$, then $\frac{1}{\tau_1}b(\tau_1{\bm \nu}) \le \frac{1}{\tau_2}b(\tau_2{\bm \nu}) \le 1,\ \forall \tau_1 \ge \tau_2 > 0,\, \forall {\bm \nu}\in [0,\infty)^N$, which gives 
\begin{align*}
\frac{1}{\tau_1}\mathcal{P}_{\tau_1}\supset \frac{1}{\tau_2}\mathcal{P}_{\tau_2},\ \forall \tau_1\ge \tau_2>0.
\end{align*} 

Since $\mathcal{P}=\bigcup_{\tau\in\mathbb{R}^+}\left(\frac{1}{\tau}\!\mathcal{P}_\tau\right)$, it is trivial that $|\mathcal{P}|>0$. Furthermore, for any ${\bm \nu} \in \mathcal{P}$, $b(\tau{\bm \nu})\le \tau$ with $\tau$ large enough. Combining with Assumption \ref{sec4:assump1}.3 that $b(\tau{\bm \nu}) > c\tau|{\bm\nu}|$, this implies $|{\bm \nu}|< \frac{1}{c}$. Thus, $\mathcal{P}\subset B({\bm 0},1/c)$ and $|\mathcal{P}| < \infty$.

If $\mathcal{P}$ Jordan measurable, since the family $(\mathcal{P}_\tau)_{\tau \in \mathbb{R}^+}$ has been proved to satisfy the conditions of Lemma \ref{point_count2}, we can apply this to get \eqref{eq:6c}.
\fproof

As seen in the proof, the continuity of $b$ (Assumption \ref{sec4:assump1}.1) assures that the superlevel sets $\mathcal{P}_{\tau}$ are ``well-behaved'' (Jordan measurable). Meanwhile, the monotonicity of $H_{\bm \nu}$ (Assumption \ref{sec4:assump1}.2) leads to the ascending (or descending) property of the chain $\left(\frac{1}{\tau}\!\mathcal{P}_\tau\right)_{\tau \in \mathbb{R}^+}$. To guarantee the limiting set $\mathcal{P}$ is bounded and not null, we assume $c|{\bm \nu}| < b({\bm \nu}) < C|{\bm \nu}|$ for some $0< c<C$ (Assumption \ref{sec4:assump1}.3), so that $B({\bm 0},1/C)\subset \mathcal{P}\subset B({\bm 0},1/c)$. It should be noted that $c$ and $C$ are generic constants, which are only utilized to represent the boundedness of $\mathcal{P}$ and do not affect our convergence rate, thus a specification of $c$ and $C$ is not necessary. In the subsequent analysis, we applies \eqref{eq:6c} to derive an error estimate, only depending on $\mathcal{P}$ and the parameter dimension $N$, of the form $M\exp(-(M/|{\mathcal{P}}|)^{1/N})$. This rate is consistent with the proven sub-exponential convergence $M \exp(-(\kappa M)^{1/N})$ for some simple coefficient upper bounds \cite{BNTT14}. Nevertheless, our analysis completely exploits { detailed information on the size and shape of the index sets} in the asymptotic regime via the introduction of $\mathcal{P}$ and, as a result, acquires the optimal value of $\kappa$. 


It is worth remarking that Lemma \ref{point_count2} requires $\mathcal{P}$ to be Jordan measurable. Indeed, we show here a simple counterexample in which $\mathcal{P}$ is not Jordan measurable and \eqref{eq:6c} fails to hold. Consider the integer-indexed collection of Jordan measurable sets $(\mathcal{P}_j)_{j\in \mathbb{N}}$ defined by 
\begin{align*}
\mathcal{P}_j = j\left([0,1]\setminus \left\{ {p}/{q}: p,q\in \mathbb{Z},\, 0\le p\le q\le j \right\}\right).
\end{align*}
Observing that $(\frac{1}{j}\mathcal{P}_j)_{j\in \mathbb{N}}$ is descending towards $\mathcal{P} = [0,1]\setminus \mathbb{Q}$, which is not Jordan measurable. We have $\# (\mathcal{P}_j\cap \mathbb{Z}^N) = 0\  \, \forall j $ while $|\mathcal{P}| = 1$, contradictory to \eqref{eq:6c}. The conditions on Jordan measurability of $\mathcal{P}$ is, however, not restrictive in the context of quasi-optimal methods, since the shapes of limiting sets are often not very {  complicated, e.g., fractal.} Indeed, all examples investigated herein show the convexity of $\mathcal{P}$, which trivially implies its Jordan measurability, as required.

The mathematical evidence that Assumption \ref{sec4:assump1} is satisfied by published Taylor and Legendre coefficient estimates will be presented in Section \ref{section:conv_anal}. Four examples of upper bounds $e^{-b({\bm \nu})}$ will be considered, including ${\bm \rho}^{-{\bm \nu}}$ (as in \eqref{taylor_est}), $ \inf_{(\delta,{\bm \rho})\in \bm{Ad}} (\frac{ {\bm \rho}^{-{\bm \nu}}}{\delta})$ (as in \eqref{bound_opt}), $ {\bm \rho}^{-{\bm \nu}}\prod_{i=1}^N\sqrt{2\nu_i+1}$ (as in \eqref{legendre_est}), and $\frac{|{\bm \nu}|!}{\bm{\nu}!}{\bm \alpha}^{{\bm \nu}} $ (as in \cite{BTNT12,CDS10}). For now, with Lemma \ref{point_count2} giving an approximation for $\# (\mathcal{P}_j\cap \mathbb{Z}^N)$, it remains to study the estimation problem \eqref{sec4:intro2}. We proceed to prove the following supporting result. 

\begin{lemma}
\label{lemma:estim}
For any $N,J,L \in \mathbb{N}$, if $J\ge  \max\left\{\frac{1}{e^{1/N}-1},\frac{L}{e^{(L-1)/N}-1}\right\}$, it gives
\begin{align}
\label{lemma:eq}
\sum_{j\ge J} j^N e^{-j} \le L J^N e^{-J}\frac{e}{e-1}.
\end{align}
Particularly, 
\begin{align}
\sum_{j\ge J} j^N e^{-j} &\le 2 J^N e^{-J}\frac{e}{e-1},\ \ \forall  J\ge \frac{2}{e^{1/N}-1}, \label{lemma:eq2}
\\
\sum_{j\ge J} j^N e^{-j} &\le (N+1) J^N e^{-J}\frac{e}{e-1},\ \ \forall  J\ge \frac{1}{e^{1/N}-1}, N\ge 4. \label{lemma:eq3}
\end{align}
\end{lemma}





\sproof
We have 
\begin{align}
\label{eq:1}
\frac{1}{J^N}  \sum_{j\ge J} j^N e^{-j} =\, \sum_{k\ge 0} \sum_{\ell=0}^{L-1} \left [ \left( 1+ \frac{Lk + \ell}{J} \right)^N e^{-J-Lk -\ell} \right ]. 
\end{align}
We prove that for every $k\ge 0$, $0\le \ell \le L-1$,
\begin{align}
\label{eq:2}
 &\ \  \left( 1+ \frac{Lk + \ell}{J} \right)^N e^{-J-Lk -\ell}  \le e^{-J-k}.
 \end{align}
 Consider $\ell = 0$. If $k=0$, \eqref{eq:2} holds trivially. If $k>0$, it is equivalent to 
 \begin{align*}
 \left( 1+ \frac{Lk}{J} \right)^N \le e^{(L-1)k},\ \ \mbox{ or} \ \ \ J\ge\frac{Lk}{e^{(L-1)k/N}-1},
 \end{align*}
which is true since $J\ge  \frac{L}{e^{(L-1)/N}-1}$.

Now, for $\ell>0$, observe that 
\begin{align*}
\left( 1+ \frac{Lk + \ell}{J} \right)^N  &=  \left( 1+ \frac{Lk}{J} \right)^N  \left( 1+ \frac{\ell}{Lk+J} \right)^N\\ 
&\le \, e^{(L-1)k}\left( 1+ \frac{\ell}{J} \right)^N \le e^{(L-1)k+\ell},  
\end{align*}
since $J\ge  \frac{1}{e^{1/N}-1}$ and \eqref{eq:2} follows. 

Combining \eqref{eq:1} and \eqref{eq:2} gives 
\begin{align*}
\frac{1}{J^N}  \sum_{j\ge J} j^N e^{-j} \le L \sum_{j\ge J} e^{-j} = L e^{-J}\frac{e}{e-1},
\end{align*}
which yields \eqref{lemma:eq}. 

\eqref{lemma:eq2} can be obtained from \eqref{lemma:eq} with $L =2 $. For \eqref{lemma:eq3}, applying \eqref{lemma:eq} with $L = N+ 1$, we only need to verify $\frac{1}{e^{1/N}-1} \ge \frac{N+1}{e-1}$. We have 
\begin{align*}
\frac{e-1}{e^{1/N}-1} = \sum_{i = 0}^{N-1} e^{i/N} \ge N+1, 
\end{align*} 
since $e^{(N-1)/N}\ge 1.5$ and $e^{(N-2)/N}\ge 1.5$ for $N\ge 4$, and $e^{i/N}\ge 1$ for all $0\le i\le N-3$, proving \eqref{lemma:eq3}. 
\fproof

It is easy to see that $\sum_{j\ge J} j^N e^{-j}$ is also bounded from below by 
\begin{align*}
\sum_{j\ge J} j^N e^{-j} \ge  J^N \sum_{j\ge J} e^{-j} \ge J^N e^{-J} \frac{e}{e-1},\ \ \forall J\in \mathbb{N},
\end{align*}
verifying the sharpness of estimate \eqref{lemma:eq}. This sub-exponential convergence rate, however, is effective with $J\ge \frac{1}{e^{1/N}-1}\simeq N$. Since $J^N e^{-J}$ is increasing with respect to $J$ for $J<N$, this seems not an appropriate rate to describe the decay of $\sum_{j\ge J} j^N e^{-j}$ in the pre-asymptotic regime. 

We are now ready to analyze the asymptotic truncation error of the general multi-indexed series $\sum_{{\bm \nu} \in \mathcal{S}} e^{-b({\bm \nu})}$ relevant to quasi-optimal {  and best $M$-term} Taylor and Legendre approximations. The main result of this section is stated and proved below. 

\begin{theorem}
\label{conv_rate1}
Consider the multi-indexed series $\sum_{{\bm \nu}\in \mathcal{S}} e^{-b({\bm \nu})}$ with $b: [0,\infty)^N\to\mathbb{R}$ satisfying Assumption \ref{sec4:assump1}. For $\tau\in (0,\infty)$, denote $\mathcal{P}_\tau = \left\{\bm{\nu} \in [0,\infty)^N:\, b(\bm{\nu})\le \tau \right\}$ and $\Lambda_M$ the set of indices corresponding to $M$ largest $e^{-b(\bm{\nu})}$. Define $\mathcal{P}=\bigcap_{\tau\in\mathbb{R}^+}\left(\frac{1}{\tau}\!\mathcal{P}_\tau\right)$ or $\mathcal{P}=\bigcup_{\tau\in\mathbb{R}^+}\left(\frac{1}{\tau}\!\mathcal{P}_\tau\right)$ as in \eqref{define.P}. If $\mathcal{P}$ is Jordan measurable, for any $\varepsilon>0$, there exists $M_{\varepsilon} >0$ depending on $\varepsilon$ such that 

\begin{align}
\label{conv_rate1:eq}
\sum\limits_{{\bm \nu}\notin \Lambda_M} e^{-b({\bm \nu})} \le C_{u}(\varepsilon)\,  M \exp \left( - \left({\frac{M}{|\mathcal{P}|(1+\varepsilon)}}\right)^{1/N} \right)
\end{align}
for all $M> M_{\varepsilon}$. Here, $C_{u}(\varepsilon) =  \left(4 e+4 \varepsilon e- {2} \right )  \frac{e}{e-1}$.
\end{theorem}

\sproof
We apply Lemma \ref{lemma:prop.of.b} to get 
\begin{align*}
|\mathcal{P}| = \lim_{\tau\to \infty} \frac{1}{\tau^{N}}\cdot \#(\mathcal{P}_\tau\cap \mathbb{Z}^N).
\end{align*}

For a fixed $\varepsilon >0$, there exists $\Delta_\varepsilon>0$ such that for all integer $ j>\Delta_\varepsilon$, 
\begin{gather}
\label{eq:3b}
\begin{aligned}
&\ \ \ -\varepsilon|\mathcal{P}| \le |\mathcal{P}| - \frac{1}{j^{N}}\cdot \# (\mathcal{P}_j \cap \mathbb{Z}^N) \le \frac{1}{2} |\mathcal{P}|, \\
\mbox{i.e.,}\ \  &\frac{1}{2} j^{N} |\mathcal{P}| \le  \# (\mathcal{P}_j \cap \mathbb{Z}^N)  \le j^{N} |\mathcal{P}| (1 +\varepsilon).
\end{aligned}
\end{gather}

To analyze the asymptotic convergence of $ \sum\limits_{{\bm \nu}\notin\Lambda_M}e^{-b({\bm \nu})}$, it is sufficient to consider this sum with $\Lambda_M = \mathcal{P}_J \cap \mathbb{Z}^N$, $J\in \mathbb{N}$. First, observe that for all integer $J> \Delta_\varepsilon$ and $J \ge \frac{1}{e^{1/N} -1 }$, from \eqref{eq:3b},
{\allowdisplaybreaks
\begin{align*}
&\sum_{{\bm \nu}\notin \mathcal{P}_J \cap \mathbb{Z}^N }e^{-b({\bm \nu})} \le \sum_{j \ge  J} ( \# (\mathcal{P}_{j+1} \cap \mathbb{Z}^N) -  \# (\mathcal{P}_{j} \cap \mathbb{Z}^N))e^{-j} 
\\
\le\, & \sum_{j \ge  J} \left [ {(j+1)}^{N} |\mathcal{P}|(1 +\varepsilon) - \frac{1}{2} j^{N} |\mathcal{P}|\right ] e^{-j}
\\
\le\, &|\mathcal{P}| \sum_{j\ge J}\left[ {(j+1)}^{N} - {j}^{N} \right] e^{-j} +   |\mathcal{P}| \sum_{j\ge J}\left[ \varepsilon {(j+1)}^{N} + \frac{1}{2}{j}^{N} \right] e^{-j}
\\
\le\, & (e-1) |\mathcal{P}|  \sum_{j \ge  J} j^{N} e^{-j} + \left (\varepsilon e+\frac{1}{2}\right)  |\mathcal{P}|  \sum_{j \ge  J} j^{N} e^{-j},
\end{align*} 
}
the last estimate coming from $(j+1)^{N} < e j^{N}$.

Apply Lemma \ref{lemma:estim} with $L = 2$ and $J\ge  \frac{2}{e^{1/N}-1}$, we have 
\begin{align}
\sum_{{\bm \nu}\notin \mathcal{P}_J \cap \mathbb{Z}^N }e^{-b({\bm \nu})} \le \, &  \left(2 e+2 \varepsilon e- {1} \right ) |\mathcal{P}|     {J}^{N} e^{-J}  \frac{e}{e-1}.  
 \label{eq:4}
\end{align}

Now, we need to write \eqref{eq:4} in term of $M = \# \Lambda_M = \# (\mathcal{P}_J\cap \mathbb{Z}^N)$. From \eqref{eq:3b}, it is easy to see that 
\begin{align}
\label{eq:5}
\frac{1}{2}  J^{N} |\mathcal{P}|  \le M \le  J^{N} |\mathcal{P}| (1+\varepsilon).
\end{align}

Combining \eqref{eq:4}--\eqref{eq:5} gives 
\begin{align*}
\sum_{{\bm \nu}\notin \Lambda_M }e^{-b({\bm \nu})} &\le \, C_u(\varepsilon) M \exp \left( - \left({\frac{M}{|\mathcal{P}|(1+\varepsilon)}}\right)^{1/N} \right),
\end{align*}
where $C_{u}(\varepsilon) =  \left(4 e+4 \varepsilon e- {2} \right )  \frac{e}{e-1}$, as desired.
\fproof

\begin{remark}[Theoretical minimum cardinality $M_{\varepsilon}$]
The error estimate \eqref{conv_rate1:eq} holds with 
\begin{align}
M> M_{\varepsilon} :=  \# (\mathcal{P}_{J_{\varepsilon}}\cap \mathbb{Z}^N),\mbox{ where }J_{\varepsilon} = \max\left\{ \frac{2}{e^{1/N}-1}, \Delta_{\varepsilon}\right\}. 
\label{J_cond}
\end{align}
It is shown in \eqref{eq:3b} that $\Delta_{\varepsilon}$ is decreasing with respect to $\varepsilon$. Thus, a stronger convergence rate, corresponding to smaller $\varepsilon$, would be realized at larger cardinality $M$. An evaluation of $\Delta_{\varepsilon}$ is not accessible to us in general, making explicit computation (or mathematical formula) of minimum cardinality $M_{\varepsilon}$ not feasible. However, in the settings where $\mathcal{P}$ is a rational convex polytope, $\Delta_{\varepsilon}$ can be acquired computationally. The interplay between $\varepsilon$ and $M_{\varepsilon}$ will be investigated through several examples within such settings in Section \ref{sec:validity}. 

In any case, \eqref{conv_rate1:eq} requires $J \ge \frac{2}{e^{1/N}-1}$. This condition can be relaxed with a slightly weaker estimate. Indeed, applying \eqref{lemma:eq3} instead of \eqref{lemma:eq2} in the proof of Theorem \ref{conv_rate1}, one gets 
\begin{align}
\label{eq_new:1}
\sum\limits_{{\bm \nu}\notin \Lambda_M} e^{-b({\bm \nu})} \le \frac{N+1}{2}C_{u}(\varepsilon)\,  M \exp \left( - \left({\frac{M}{|\mathcal{P}|(1+\varepsilon)}}\right)^{1/N} \right),
\end{align} 
given
\begin{align} 
\label{J'_cond}
M> M'_{\varepsilon} := \# (\mathcal{P}_{{J}'_{\varepsilon}}\cap \mathbb{Z}^N),\mbox{ where }{J}'_{\varepsilon} = \max\left\{ \frac{1}{e^{1/N}-1}, \Delta_{\varepsilon}\right\}.
\end{align} 
\end{remark}

\begin{remark}[An extension of Theorem \ref{conv_rate1} {  for $b({\bm \nu}) \in \Theta(|{\bm \nu}|^{\beta})$}] 
The convergence estimate \eqref{conv_rate1:eq} does not apply for $|\mathcal{P}| = 0$ or $\mathcal{P}$ unbounded ($b({\bm \nu}) \notin \Theta(|{\bm \nu}|)$). With minor modifications in the above analysis, our results can be extended to a wider class of $b({\bm \nu})$ where Assumption \ref{sec4:assump1}.3 ($b({\bm \nu}) \in \Theta(|{\bm \nu}|)$) is replaced by the condition that $b({\bm \nu}) \in \Theta(|{\bm \nu}|^{\beta})$ (i.e., there exist constants $0<c<C$ such that $c|{\bm \nu}|^{\beta} < b({\bm \nu}) < C|{\bm \nu}|^{\beta}$ as ${\bm \nu} \to \infty$) with some fixed $\beta>0$. In such cases, it gives 
\begin{align*}
\sum\limits_{{\bm \nu}\notin \Lambda_M} e^{-b({\bm \nu})} \le C_{u}(\varepsilon)\,  M \exp \left( - \left({\frac{M}{|\mathcal{P}|(1+\varepsilon)}}\right)^{\beta/N} \right)
\end{align*} 
as $M\to\infty$. Here, $\mathcal{P}=\bigcap_{\tau\in\mathbb{R}^+}\left(\frac{1}{\tau^{1/\beta}}\!\mathcal{P}_\tau\right)$ or $\mathcal{P}=\bigcup_{\tau\in\mathbb{R}^+}\left(\frac{1}{\tau^{1/\beta}}\!\mathcal{P}_\tau\right)$ (depending on whether $\frac{1}{\tau^{1/\beta}}\!\mathcal{P}_\tau$ is descending or ascending). 
\end{remark} 


\section{Asymptotic convergence rates of quasi-optimal { and best $M$-term} approximations} 
\label{section:conv_anal}
As we have seen so far, the error of a quasi-optimal polynomial approximation can be estimated by {  the series of corresponding coefficient upper bounds.} {  This also represents an accessible convergence estimate for the best $M$-term approximation, as discussed in Section \ref{sec:intro}. We will verify in this section that for most upper bounds developed in recent publications, such series fall into the class of multi-indexed series analyzed in Section \ref{anal_general_series}. Particularly, in all considered cases, the coefficient estimates, written as $e^{-b({\bm \nu})}$, will be proved to satisfy Assumption \ref{sec4:assump1} and $\sum_{{\bm \nu}\in \mathcal{S}} e^{-b({\bm \nu})}$ can be treated by Theorem \ref{conv_rate1}.} 

Given a vector ${\bm \rho} = (\rho_i)_{i\le 1\le N}$ with $\rho_i > 1\ \forall i$, we define ${\bm \lambda} = (\lambda_i)_{1\le i \le N}$ such that 
\begin{align*}
\lambda_i = \log \rho_i > 0 \  \, \forall\, 1\le i \le N. 
\end{align*}
{  In Section \ref{sec:quasi_taylor} and \ref{sec:quasi_legendre}, we study the error analysis of quasi-optimal and best $M$-term Taylor and Legendre approximations. (However, for ease of presentation, in what follows, we mostly refer to the analysis as error estimate of quasi-optimal methods).} A computational comparison of our proposed estimate with existing results is showing in Section \ref{sec:comparison}. 

\subsection{Error analysis of quasi-optimal Taylor approximations} 
\label{sec:quasi_taylor}
We start with the quasi-optimal methods corresponding to a basic coefficient bound of the form ${\bm \rho}^{-{\bm \nu}}$ (see Proposition \ref{theorem:taylor_est}). These are reasonable schemes for Taylor approximations of elliptic problems with the random fields composed of non-overlapping basis functions. The convergence result is stated in the following proposition.




\begin{proposition}
\label{theorem:aniso}
Consider the Taylor series $\sum\limits_{{\bm \nu}\in \mathcal{S}} t_{\bm \nu}{\bm y}^{\bm \nu}$ of $u$. Assume that
\begin{align}
\label{taylor_est2}
\|t_{\bm \nu}\|_{V(D)} \le \frac{\|f\|_{V^*(D)}}{\delta}{\bm \rho}^{-{\bm \nu}} 
\end{align}
holds for all ${\bm \nu}\in\mathcal{S}$, as in Proposition \ref{theorem:taylor_est}. Denote by $\Lambda_M$ the set of indices corresponding to $M$ largest bounds in \eqref{taylor_est2}. For any $\varepsilon > 0$, there exists $M_{\varepsilon} > 0$ depending on $\varepsilon$ such that
\begin{align}
\label{theorem:est1}
\sup_{{\bm y}\in\Gamma} \left\| u({\bm y}) -\!\! \!\sum\limits_{{\bm \nu} \in \Lambda_M}  t_{\bm \nu}{\bm y}^{\bm \nu} \right\|_{V(D)} \le \frac{\|f\|_{V^*(D)}}{\delta} C_u(\varepsilon) M \exp\! \left(\! -\!\left({\frac{MN! \prod_{i=1}^{N}\lambda_i}{(1+\varepsilon)}}\right)^\frac{1}{N} \right)
\end{align}
 for all $M>M_{\varepsilon}$. 
\end{proposition}

\sproof
We have by triangle inequality
\begin{align}
\sup_{{\bm y}\in\Gamma} \left\| u({\bm y}) - \sum\limits_{{\bm \nu} \in \Lambda_M}  t_{\bm \nu}{\bm y}^{\bm \nu} \right\|_{V(D)} \le \sum_{{\bm \nu}\notin \Lambda_M} \|t_{\bm \nu}\|_{V(D)} \le   \frac{\|f\|_{V^*(D)}}{\delta} \sum_{{\bm \nu}\notin \Lambda_M} {\bm \rho}^{-{\bm \nu}}.  \label{triag_ineq}
\end{align}

For $ {\bm \nu} \in[0,\infty)^N$, define $b({\bm \nu}) = \sum\limits_{i=1}^N \lambda_i \nu_i$, so that $  {\bm \rho}^{-{\bm \nu}} = e^{-b({\bm \nu})}\ \, \forall  {\bm \nu} \in \mathcal{S}$. We notice that the quasi-optimal index sets in this case are the Total Degree spaces:
\begin{align*}
\mathcal{P}_j\cap \mathbb{Z}^N = \left\{{\bm \nu} \in \mathcal{S}: {\bm \rho}^{-{\bm \nu}} \ge e^{-j} \right\} = \left\{{\bm \nu} \in \mathcal{S}: \sum\limits_{i=1}^N \lambda_i\nu_i \le j \right\},\ \forall j \in \mathbb{N}.
\end{align*}

Since $\lambda_i>0\ \, \forall i$, it is easy to check that the map $b$ satisfies Assumption \ref{sec4:assump1} with $H_{\bm \nu}$ being constant $\forall {\bm \nu}$. Observing that $\mathcal{P} = \bigcap_{\tau\in\mathbb{R}^+}\left(\frac{1}{\tau}\mathcal{P}_\tau\right) =  \{{\bm \nu} \in [0,\infty)^N: \sum_{i=1}^N \lambda_i\nu_i \le 1 \}$, we can specify $|\mathcal{P}| = \frac{1}{N!(\lambda_1\ldots \lambda_N)}$. 

We are now ready to apply Theorem \ref{conv_rate1} to obtain
\begin{align*}
\sum_{{\bm \nu}\notin \Lambda_M }{\bm \rho}^{-{\bm \nu}} &\le \, C_u(\varepsilon) M \exp \left( - \left({\frac{M}{|\mathcal{P}|(1+\varepsilon)}}\right)^{1/N} \right)
\\
&\le \, C_u(\varepsilon) M \exp \left( -\left({\frac{MN! \prod_{i=1}^{N}\lambda_i}{(1+\varepsilon)}}\right)^{1/N} \right),
\end{align*}
which proves \eqref{theorem:est1}.
\fproof

We proceed to analyze the quasi-optimal Taylor approximations based on best analytical bound provided by Proposition \ref{theorem:taylor_est}. Although this method is not easily implementable, an asymptotic error estimate can be obtained as a simple corollary of Theorem \ref{conv_rate1}. It is reasonable to assume that the set $\textit{\textbf{Ad}}$ of all admissible $(\delta,{\bm \rho})$ is bounded: as seen through several examples in Figure \ref{domUE}, the domains of uniform ellipticity do not expand infinitely in complex plane.  

\begin{proposition}
\label{best_bound}
Consider the Taylor series $\sum\limits_{{\bm \nu}\in \mathcal{S}} t_{\bm \nu}{\bm y}^{\bm \nu}$ of $u$. Assume
\begin{align}
\label{theorem:est3}
\|t_{\bm \nu}\|_{V(D)} \le \inf\limits_{(\delta,{\bm \rho})\in \bm{Ad}} \frac{\|f\|_{V^*(D)}}{\delta} {\bm \rho}^{-{\bm \nu}}
\end{align}
holds for all ${\bm \nu}\in\mathcal{S}$, as in \eqref{bound_opt}, with $\textbf{Ad}$ being bounded. Denote by $\Lambda_M$ the set of indices corresponding to $M$ largest bounds in \eqref{theorem:est3}. For any $\varepsilon>0$, there exists $M_{\varepsilon}>0$ depending on $\varepsilon$ such that 
\begin{align}
\label{theorem:est4}
\sup_{{\bm y}\in\Gamma} \left\| u({\bm y}) - \sum\limits_{{\bm \nu} \in \Lambda_M}  t_{\bm \nu}{\bm y}^{\bm \nu} \right\|_{V(D)}\!\le   {\|f\|_{V^*(D)}} C_u(\varepsilon) M \exp\! \left( - \left({\frac{M}{|\mathcal{P}|(1+\varepsilon)}}\right)^{1/N} \right)
\end{align}
for all $M > M_{\varepsilon} $. Here, 
$\mathcal{P} = \left\{{\bm \nu} \in [0,\infty)^N:  \sum_{i=1}^N (\log \rho_i) \nu_i \le 1 \ \ \forall (\delta,{\bm \rho})\in \bm{Ad}\right\}.$
\end{proposition}

\sproof
First, we have 
\begin{align}
\sup_{{\bm y}\in\Gamma} \left\| u({\bm y}) - \sum\limits_{{\bm \nu} \in \Lambda_M}  t_{\bm \nu}{\bm y}^{\bm \nu} \right\|_{V(D)} \le \sum_{{\bm \nu}\notin \Lambda_M} \|t_{\bm \nu}\|_{V(D)} \le   {\|f\|_{V^*(D)}}\sum_{{\bm \nu}\notin \Lambda_M}  \inf\limits_{(\delta,{\bm \rho})\in \bm{Ad}} \frac{{\bm \rho}^{-{\bm \nu}}}{\delta} .  \label{triag_ineq2}
\end{align}

Recall that $\lambda_i = \log \rho_i\ \forall i$. With abuse of notation, we say $(\delta,{\bm \lambda}) \in \bm{Ad}$ iff $(\delta,{\bm \rho}) \in \bm{Ad}$. For $ {\bm \nu} \in[0,\infty)^N$, define $b({\bm \nu}) = \sup\limits_{(\delta,{\bm \lambda})\in \bm{Ad}}\left(\log\delta + \sum\limits_{i=1}^N \lambda_i \nu_i\right)$, so that $  \inf\limits_{(\delta,{\bm \rho})\in \bm{Ad}} \frac{{\bm \rho}^{-{\bm \nu}}}{\delta} = e^{-b({\bm \nu})}\ \, \forall  {\bm \nu} \in \mathcal{S}$. The quasi-optimal index sets in this case are:
\begin{align*}
\mathcal{P}_j\cap \mathbb{Z}^N =  \left\{{\bm \nu} \in \mathcal{S}: \sup\limits_{(\delta,{\bm \lambda})\in \bm{Ad}}\left(\log\delta + \sum\limits_{i=1}^N \lambda_i \nu_i\right) \le j \right\},\ \forall j\in \mathbb{N}.
\end{align*}

We will show that $b$ fulfills Assumption \ref{sec4:assump1}. It is easy to check that $b$ is convex. As a consequence, for any $\tau>0,\, \mathcal{P}_\tau$ and $\bigcap_{\tau\in\mathbb{R}^+} \left(\frac{1}{\tau} \mathcal{P}_\tau\right)$ are convex (and Jordan measurable). Since $\textit{\textbf{Ad}}$ is bounded, there exist $0<c<C$ such that $c|{\bm \nu}| <b({\bm \nu}) < C|{\bm \nu}|$ as ${\bm \nu} \to \infty$. Now, let $\tau\ge \tau' >0$, it gives
\begin{align*}
\frac{1}{\tau'}\left(\log\delta + \sum\limits_{i=1}^N \lambda_i \tau'\nu_i\right) \le \frac{1}{\tau}\left(\log\delta + \sum\limits_{i=1}^N \lambda_i \tau\nu_i\right),\ \forall (\delta,{\bm \lambda})\in \bm {Ad},\, {\bm \nu}\in [0,\infty)^N,
\end{align*}
since $\delta < 1$. Hence, $H_{\bm \nu}({\tau'}) \le H_{\bm \nu}({\tau}),\, \forall {\bm \nu}\in [0,\infty)^N $. 

We can apply Theorem \ref{conv_rate1} to get the asymptotic estimate
\begin{align}
\sum_{{\bm \nu}\notin \Lambda_M}  \inf\limits_{(\delta,{\bm \rho})\in \bm{Ad}} \frac{{\bm \rho}^{-{\bm \nu}}}{\delta} \le C_u(\varepsilon) M \exp \left( - \left({\frac{M}{|\mathcal{P}|(1+\varepsilon)}} \right)^{1/N} \right), \label{eq:9}
\end{align} 
where $\mathcal{P}:=  \bigcap_{\tau\in\mathbb{R}^+} \left(\frac{1}{\tau} \mathcal{P}_\tau\right)= \left\{{\bm \nu} \in [0,\infty)^N:  \sum_{i=1}^N \lambda_i \nu_i \le 1 \ \ \forall (\delta,{\bm \lambda})\in \bm{Ad}\right\}$. Combining \eqref{triag_ineq2} and \eqref{eq:9} gives \eqref{theorem:est4}, concluding the proof. 
\fproof


\subsection{Error analysis of quasi-optimal Legendre approximations}

\label{sec:quasi_legendre}

For the first example, we consider quasi-optimal methods for Legendre approximations of elliptic PDEs with the random field consisting of basis functions with disjoint supports. In \cite{BNTT14}, these problems were computationally treated with bounds of type \eqref{taylor_est2} and Total Degree index sets with some success. However, those bounds are not analytically optimal, as the true exponential decay of coefficients is penalized by a large multiplier. In the following, we establish a convergence analysis for the sharper upper bound $ {\bm \rho}^{-{\bm \nu}}\prod_{i=1}^N\sqrt{2\nu_i+1}$ of Legendre coefficients (see Section \ref{estim_coef}). Whether the quasi-optimal method corresponding to this estimate outperforms Total Degree approximations in computation is an interesting subject to study next.  

\begin{proposition}
\label{theorem:legendre}
Consider the Legendre series $\sum\limits_{{\bm \nu}\in \mathcal{S}} v_{\bm \nu}L_{\bm\nu}$ of $u$. Assume that
\begin{align}
\label{legendre_est2}
\|v_{\bm \nu}\|_{V(D)} \le C_{{\bm \rho},\delta} {\bm \rho}^{-{\bm \nu}}\prod_{i=1}^N\sqrt{2\nu_i+1}
\end{align}
holds for all ${\bm \nu}\in\mathcal{S}$, as in Proposition \ref{theorem:legendre_est}. Denote by $\Lambda_M$ the set of indices corresponding to $M$ largest bounds in \eqref{legendre_est2}. For any $\varepsilon > 0$, there exists a constant $M_{\varepsilon} > 0$ depending on $\varepsilon$ such that 
\begin{align}
\label{theorem:est2}
 \left\| u - \sum\limits_{{\bm \nu} \in \Lambda_M}  v_{\bm \nu}L_{\bm \nu} \right\|_{V(D)\otimes L^2_{\varrho}(\Gamma)}^2\!\! \le C_{{\bm \rho},\delta}^2 C_u(\varepsilon)    M \exp\! \left(\! - 2 \left({\frac{MN! \prod_{i=1}^{N}\lambda_i}{(1+\varepsilon)}}\right)^{1/N} \right)
\end{align}
for all $M>M_{\varepsilon}$. 
\end{proposition}

\sproof
First, we have 
\begin{align*}
\left\|u - \sum_{{\bm \nu} \in \Lambda_M} v_{\bm \nu}L_{\bm \nu}\right\|_{V(D)\otimes L^2_{\varrho}(\Gamma)}^2  = \ \sum_{{\bm \nu} \notin \Lambda_M} {\left \| v_{\bm \nu}\right \|^2_{V(D)} \le \ C_{{\bm \rho},\delta}^2 \sum_{{\bm \nu} \notin \Lambda_M}  {\bm \rho}^{-2{\bm \nu}}\prod_{i=1}^N({2\nu_i+1}) }.
\end{align*}
For $ {\bm \nu} \in[0,\infty)^N$, define $b({\bm \nu}) = \sum\limits_{i=1}^N \left( 2 \lambda_i \nu_i - \log(2\nu_i + 1) \right)$, so that \\$  {\bm \rho}^{-2{\bm \nu}}\prod_{i=1}^N({2\nu_i+1})$ $ = e^{-b({\bm \nu})}\ \, \forall  {\bm \nu} \in \mathcal{S}$. We notice that the quasi-optimal index sets in this case given by:
\begin{align*}
\mathcal{P}_j\cap \mathbb{Z}^N  = \left\{{\bm \nu} \in \mathcal{S}:  \sum\limits_{i=1}^N \left( 2 \lambda_i \nu_i - \log(2\nu_i + 1) \right) \le j \right\},\ \forall j\in \mathbb{N}.
\end{align*}

We proceed to prove $b$ satisfies Assumption \ref{sec4:assump1}. It is easy to check that $b({\bm \nu})$ is continuous. As ${\bm \nu}\to \infty$, $\lambda_{\min}|{\bm \nu}| < b({\bm \nu}) < 2 \lambda_{\max}|{\bm \nu}|$, where $\lambda_{\min} = \min_{1\le i\le N} \lambda_i$ and $\lambda_{\max} = \max_{1\le i\le N} \lambda_i$. 
Also, observing $\log(at+ 1) \ge t\log(a+1)$ for every  $a\ge 0,\, 0\le t \le 1$, we have 
\begin{align*}
H_{\bm \nu}(\tau') &=\! \sum\limits_{i=1}^N \left( 2 \lambda_i \nu_i \!-\! \frac{1}{\tau'} \log(2\tau'\nu_i + 1) \right) \\
&\le \sum\limits_{i=1}^N \left( 2 \lambda_i \nu_i\! -\! \frac{1}{\tau} \log(2\tau\nu_i + 1) \right) 
=\! H_{\bm \nu}(\tau)
\end{align*} 
for all $\tau,\tau'\in (0,\infty),\, \tau\ge \tau'$. 

{  It is easy to see that $\mathcal{P} := \bigcap\limits_{\tau\in\mathbb{R}^+}\left(\frac{1}{\tau}\mathcal{P}_\tau\right) = \left\{{\bm \nu} \in [0,\infty)^N: \sum\limits_{i=1}^N 2 \lambda_i\nu_i \le 1 \right\}$. Thus, $\mathcal{P}$ is Jordan measurable and $|\mathcal{P}| = \frac{1}{2^N N!\prod_{i=1}^N \lambda_i}$. Applying Theorem \ref{conv_rate1}, we obtain  
\begin{align*}
\sum_{{\bm \nu} \notin \Lambda_M}  {\bm \rho}^{-2{\bm \nu}}\prod_{i=1}^N({2\nu_i+1}) 
 &\le \, C_u(\varepsilon) M \exp \left( - 2 \left({\frac{MN!\prod_{i=1}^N \lambda_i}{1+\varepsilon}} \right)^{1/N}\right)
\end{align*}
for all $M > M_{\varepsilon}$. This concludes our proof. }
\fproof

We remark that while the bound \eqref{legendre_est2} is weaker than \eqref{taylor_est2}, its corresponding index sets are descending towards Total Degree sets. As a result, we are able to obtain the same convergence rate as Taylor approximations. 

Now, we apply our framework to prove a convergence estimate for quasi-optimal Legendre approximations based on the coefficient exponential decay $\|v_{\bm \nu}\|_{V(D)} \le {\|f\|_{V^*(D)}}$ $\frac{|{\bm \nu}|!}{\bm{\nu}!}{\bm \alpha}^{{\bm \nu}}$. Unlike other upper bounds discussed so far, this decay is established by real analysis argument \cite{CDS10}. In the case of affine linear random fields, 
i.e.~$a(x,{\bm y}) = a_0(x) + \sum_{i=1}^N y_i \psi_i(x)$, ${\bm \alpha} = (\alpha_i)_{1\le i\le N}$ is specified by $\alpha_i = \frac{\|\psi_i\|_{L^{\infty}(D)}}{a_{\min}\sqrt{3}}$. A development and implementation of quasi-optimal method can be found in \cite{BTNT12}; however, no error estimate has been provided. In the following result, similar to the aforementioned works, we assume $\sum_{i=1}^N \alpha_i < 1$, which is necessary for the summability of sequence $\left(\frac{|{\bm \nu}|!}{\bm{\nu}!}{\bm \alpha}^{{\bm \nu}}\right)_{{\bm \nu}\in\mathcal{S}}$.


\begin{proposition}
\label{theorem:TD-FC}
Consider the Legendre series $\sum_{{\bm \nu}\in \mathcal{S}} v_{\bm \nu}L_{\bm\nu}$ of $u$. Assume there exists a vector ${\bm \alpha} = (\alpha_i)_{1\le i\le N}$ with $\alpha_i>0\ \, \forall i$ and $\sum_{i=1}^N \alpha_i < 1$ such that
\begin{align}
\label{legendre_TD-FC}
\|v_{\bm \nu}\|_{V(D)} \le {\|f\|_{V^*(D)}}\frac{|{\bm \nu}|!}{\bm{\nu}!}{\bm \alpha}^{{\bm \nu}}
\end{align}
for all ${\bm \nu}\in\mathcal{S}$. Denote by $\Lambda_M$ the set of indices corresponding to $M$ largest bounds in \eqref{legendre_TD-FC}. For any $\varepsilon>0$, there exists a constant $M_{\varepsilon} > 0$ depending on $\varepsilon$ such that 
\begin{align}
\label{eq:TD-FC}
 \left\| u - \sum\limits_{{\bm \nu} \in \Lambda_M}  v_{\bm \nu}L_{\bm \nu} \right\|_{V(D)\otimes L^2_{\varrho}(\Gamma)}^2\!\le  {\|f\|_{V^*(D)}^2}  C_u(\varepsilon)   M \exp\! \left( - \left({\frac{M}{|\mathcal{P}|(1+\varepsilon)}} \right)^{1/N}\right)
\end{align}
for all $M > M_{\varepsilon}$. Here, $\mathcal{P}  = \left\{{\bm \nu} \in (0,\infty)^N :  \sum\limits_{i=1}^N \lambda_i \nu_i - \log \frac{ |{\bm \nu}|^{|{\bm \nu}| }}{ \prod_{i=1}^N { {\nu_i}^{\nu_i}}}  < \frac{1}{2} \right\}$. 
\end{proposition}


\sproof
From \eqref{legendre_TD-FC}, we have 
\begin{align*}
\left\|u - \sum_{{\bm \nu} \in \Lambda_M} v_{\bm \nu}L_{\bm \nu}\right\|_{V(D)\otimes L^2_{\varrho}(\Gamma)}^2  = \ \sum_{{\bm \nu} \notin \Lambda_M} {\left \| v_{\bm \nu}\right \|^2_{V(D)} \le \  {\|f\|_{V^*(D)}^2}\sum_{{\bm \nu} \notin \Lambda_M}  {\bm \alpha}^{2{\bm \nu}}\left(\frac{|{\bm \nu}|!}{\bm{\nu}!}\right)^2 }. 
\end{align*}
Let $\lambda_i = - \log \alpha_i >0\ \, \forall\, 1\le  i \le N$ and $\Gamma $ denote the gamma function. Also, let $\psi_0$, $\psi_1$ and $\psi_2$ be the di-, tri- and tetra-gamma functions respectively: $\psi_0 = (\log \Gamma)',\, \psi_1 = \psi'_0 = (\log \Gamma)'',\, \psi_2 = \psi'_1 = (\log\Gamma)'''$. For $ {\bm \nu} \in[0,\infty)^N$, define $b({\bm \nu}) =  2\sum\limits_{i=1}^N \lambda_i \nu_i $ $- 2\log\frac{\Gamma(|{\bm \nu}| +  1)}{\prod_{i=1}^N{\Gamma(\nu_i + 1)}} $, so that $  {\bm \alpha}^{2{\bm \nu}}\left(\frac{|{\bm \nu}|!}{\bm{\nu}!}\right)^2 $ $ = e^{-b({\bm \nu})}\ \, \forall  {\bm \nu} \in \mathcal{S}$. The quasi-optimal index sets in this case are given by:
\begin{align*}
\mathcal{P}_j\cap \mathbb{Z}^N  = \left\{{\bm \nu} \in \mathcal{S}:  \sum\limits_{i=1}^N \lambda_i \nu_i -\log\frac{\Gamma(|{\bm \nu}| +  1)}{\prod_{i=1}^N{\Gamma(\nu_i + 1)}}  \le \frac{j}{2} \right\}.
\end{align*}

We proceed to prove $b$ satisfies Assumption \ref{sec4:assump1}. First, since $\sum_{i=1}^N\alpha_i <$ $ 1$, one can find $p\in(0,1)$ such that $\sum_{i=1}^N\alpha_i^p < 1$ and, by Theorem 7.2 in \cite{CDS10}, have $\left(\frac{|{\bm \nu}|!}{\bm{\nu}!}{\bm \alpha}^{p{\bm \nu}}\right)_{{\bm \nu}\in \mathcal{S}}$  $\ell^1$-summable. This gives   
$ \left(\frac{|{\bm \nu}|!}{\bm{\nu}!}\right)^{2} {\bm \alpha}^{2p{\bm \nu}} < 1 \mbox{ as } {\bm \nu}\to \infty$ and there follows
\begin{align*}
\left(\frac{|{\bm \nu}|!}{\bm{\nu}!}\right)^{2} {\bm \alpha}^{2{\bm \nu}} < {\bm \alpha}^{(2-2p){\bm \nu}},\mbox{ i.e.},\  b({\bm \nu}) > (2-2p)\sum_{i=1}^N \lambda_i \nu_i
\ \mbox{ as }\ {\bm \nu} \to \infty.
 \end{align*}
 
Next, define $g(\tau) = \frac{1}{\tau} \log \left(\frac{\Gamma(|{\tau{\bm \nu}}| +  1)}{\prod_{i=1}^N{\Gamma({\tau\nu_i} + 1)}} \right)$ to be a mapping from $(0,\infty)$ to $\mathbb{R}$. We will prove $H_{\bm \nu}$ is decreasing by showing $g$ is an increasing function. Observing that 
\begin{align*} 
g(\tau) = \sum\limits_{q =2}^N \frac{1}{\tau}\log\left( \frac{{\Gamma\left({\tau\sum\limits_{i=1}^{q}\nu_i} + 1\right)}}{{{\Gamma\left({\tau\sum\limits_{i=1}^{q-1}\nu_i} + 1\right)}}{{\Gamma\left({\tau\nu_q} + 1\right)}}} \right),
\end{align*}
without loss of generality, we can assume $N=2$. Consider the first derivative of $g$: 
\begin{align*}
g'(\tau) &= -\frac{1}{\tau^2} \log \left(\frac{\Gamma({\tau{\nu_1 + \tau\nu_2}} +  1)}{{\Gamma({\tau\nu_1} + 1)} {\Gamma({\tau\nu_2} + 1)} } \right) + \frac{1}{\tau^2}\frac{\Gamma ' (\tau\nu_1 + \tau\nu_2 + 1)}{\Gamma  (\tau\nu_1 + \tau\nu_2 + 1)} (\tau\nu_1 + \tau\nu_2)
\\
 &\qquad\qquad -  \frac{1}{\tau^2}\frac{\Gamma ' (\tau\nu_1  + 1)}{\Gamma  (\tau\nu_1 + 1)} \tau\nu_1 -  \frac{1}{\tau^2}\frac{\Gamma ' (\tau\nu_2  + 1)}{\Gamma  (\tau\nu_2 + 1)} \tau\nu_2 . 
\end{align*}
Then $g'(\tau)\ge 0\ \, \forall \tau >0$ iff $h( \nu_1 + \nu_2 ) \ge h(\nu_1) + h(\nu_2),\ \forall \nu_1,\nu_2\ge 0$, where $h(s) := s \psi_0(s+1) - \log(\Gamma( s+1))$. 

We have $h''(s) = s\psi_2(s+1) + \psi_1(s+1) > 0$ for any $s\ge 0$, see Theorem 1, \cite{EL00}, so $h$ is convex. Combining with the fact that $h(0) = 0$, this implies the superadditivity of $h$ in $[0,\infty)$, as desired. Note that for ${\bm \nu} \in (0,\infty)^N$, $g$ is strictly increasing in $(0,\infty)$.  
 
Since $H_{\bm \nu}$ is decreasing, define the limiting set $\mathcal{P} = \bigcup_{\tau\in \mathbb{R}^+} \left(\frac{1}{\tau}\mathcal{P}_{\tau}\right)$. We will characterize $\mathcal{P}$ and show it is Jordan measurable. Without loss of generality, we can ignore the set of points of $\mathcal{P}$ in the coordinate hyperplanes, since it is of measure zero. Using the strictly increasing property of $g$ for ${\bm \nu} \in (0,\infty)^N$, it gives
\begin{align*}
 \bigcup\limits_{\tau\in \mathbb{R}^+} \left(\frac{1}{\tau}\mathcal{P}_{\tau}\right) =  \left\{{\bm \nu} \in (0,\infty)^N :  \sum\limits_{i=1}^N \lambda_i \nu_i - \lim\limits_{\tau\to \infty}\frac{1}{\tau}\log\frac{\Gamma(|\tau{\bm \nu}| +  1)}{\prod_{i=1}^N{\Gamma(\tau\nu_i + 1)}}  < \frac{1}{2} \right\}
\end{align*}
Applying Stirling's formula, see, e.g., \cite{AS65}, yields 
\begin{align*}
 &\lim\limits_{\tau\to \infty}\frac{1}{\tau}\log\frac{\Gamma(|\tau{\bm \nu}| +  1)}{\prod_{i=1}^N{\Gamma(\tau\nu_i + 1)}} = \log \lim\limits_{\tau\to\infty} \left(  \frac{|\tau {\bm \nu}|^{|\tau {\bm \nu}| + \frac{1}{2}} e^{- |\tau {\bm \nu}|} (2\pi)^{\frac{1}{2}}}{\prod_{i=1}^N {(\tau {\nu_i})^{\tau\nu_i+\frac{1}{2}} } e^{- \tau { \nu_i}} (2\pi)^{\frac{1}{2}}}  \right)^{\frac{1}{\tau}} 
 \\
 = \, & \log \lim\limits_{\tau\to\infty}   \frac{\tau^{|{\bm \nu}| + \frac{1}{2\tau}} |{\bm \nu}|^{|{\bm \nu}| + \frac{1}{2\tau}} }{ \tau^{|{\bm \nu}| + \frac{N}{2\tau}} \prod_{i=1}^N { {\nu_i}^{\nu_i+\frac{1}{2\tau}} }  }   (2\pi)^{\frac{1-N}{2\tau} } = \, \log \frac{ |{\bm \nu}|^{|{\bm \nu}| }}{ \prod_{i=1}^N { {\nu_i}^{\nu_i}}},
\end{align*}
and we obtain 
\begin{align*}
\mathcal{P} =  \left\{{\bm \nu} \in (0,\infty)^N :  \sum\limits_{i=1}^N \lambda_i \nu_i - \log \frac{ |{\bm \nu}|^{|{\bm \nu}| }}{ \prod_{i=1}^N { {\nu_i}^{\nu_i}}}  < \frac{1}{2} \right\}.
\end{align*}

For the Jordan measurability of $\mathcal{P}$, we prove $\mathcal{P}$ is convex. It is enough to show the function $G({\bm \nu}) := \log \frac{ |{\bm \nu}|^{|{\bm \nu}| }}{ \prod_{i=1}^N { {\nu_i}^{\nu_i}}}$ is concave in $(0,\infty)^N$. Denote by $\nabla^2 G$ the Hessian matrix of $G$ and again assume $N=2$, we have 
\[ 
\nabla^2 G
= \left( \begin{array}{cc}
1/(\nu_1+\nu_2) - 1/\nu_1 & 1/(\nu_1+\nu_2) \\
1/(\nu_1+\nu_2) & 1/(\nu_1+\nu_2) - 1/\nu_2 
\end{array} \right).
\] 
Let ${\bm x} = \left( \begin{array}{c}
\! x_1 \! \\ \! x_2 \! \end{array} \right)\in \mathbb{R}^2 \setminus \{\bm 0\}$, it gives  
$ {\bm x}^\top (\nabla^2 G) {\bm  x} = \frac{(x_1 + x_2)^2}{\nu_1 + \nu_2 }  - \frac{x_1^2}{\nu_1} - \frac{x_2^2}{\nu_2 } \le 0$, by employing Cauchy-Schwarz inequality. Thus, $\nabla^2 G$ is negative semidefinite, which implies the concavity of $G$.

We can apply Theorem \ref{conv_rate1} to get the asymptotic estimate
\begin{align*}
\sum_{{\bm \nu}\notin \Lambda_M} {\bm \alpha}^{2{\bm \nu}}\left(\frac{|{\bm \nu}|!}{\bm{\nu}!}\right)^2 \le C_u(\varepsilon) M \exp \left( - \left({\frac{M}{|\mathcal{P}|(1+\varepsilon)}}\right)^{1/N} \right). 
\end{align*} 
The proof is now complete. 
\fproof 

\subsection{A computational comparison of our proposed estimate with previously established rates of convergence}
\label{sec:comparison}
Most of the established explicit error estimates for {  best $M$-term and} quasi-optimal approximations concern the coefficient bounds of the form 
\begin{align}
\|t_{\bm \nu}\|_{V(D)} \le {\bm \rho}^{-{\bm \nu}}, 
\label{compare:eq1}
\end{align} 
therefore are derived via {  the truncation error of {  $\sum_{{\bm \nu}\in \mathcal{S}}  {\bm \rho}^{-{\bm \nu}} $}}. 
We compare our approach with others in current literature in estimating this {  quantity}. Recall, we proved in Proposition \ref{theorem:aniso} that
\begin{align}
\label{compare:est1}
{  \sum_{{\bm \nu}\notin \Lambda_{M}} {\bm \rho}^{-{\bm \nu}}  } \le C_u(\varepsilon) M \exp\! \left(\! -\!\left({\frac{MN! \prod_{i=1}^{N}\lambda_i}{(1+\varepsilon)}}\right)^\frac{1}{N} \right).
\end{align} 

{{ Application of the Stechkin inequality gives 
\begin{align}
\label{CDS_bound}
\tag{stech}
\sum_{{\bm \nu}\notin \Lambda_{M}} {\bm \rho}^{-{\bm \nu}}    \le \left(\prod\limits_{i=1}^N\frac{1}{1-e^{-p\lambda_i}}\right)^{1/p} M^{1-\frac{1}{p}},
\end{align}
for every $0<p < 1$. We note that \eqref{CDS_bound} holds for every $M$ and is not asymptotic.} 

Development due to \cite{BNTT14} computes $p\in(0,1)$ minimizing \eqref{CDS_bound} for each $M$ and obtains
\begin{align}
\label{BNTT_bound}
\tag{optim}
{  \sum_{{\bm \nu}\notin \Lambda_{M}} {\bm \rho}^{-{\bm \nu}}  } \le  M \exp \left( -{\frac{1}{e}{\left(M \prod_{i=1}^{N}\lambda_i\right)^{1/N} N\xi}} \right),
\end{align}
where $\xi$ is the rate adjusting parameter varying from $0$ to $(e-1)/e$. Large $\xi$ gives stronger convergence but also require more restrictive minimum cardinality. The best convergence is only guaranteed in the limit $M\to \infty$.

Figure \ref{compareBound} shows a comparison of our error estimate with \eqref{CDS_bound} and \eqref{BNTT_bound} in computing the series $\sum_{{\bm \nu} \in \mathcal{S}} e^{-(\nu_1 + \nu_2 + 2\nu_3 + 4\nu_4)}$. \eqref{BNTT_bound} is plotted at its best possible rate with $\xi = (e-1)/e$. We also plot the exact value of $\sum_{{\bm \nu} \notin \Lambda_M} e^{-b({\bm \nu})}$, which can be calculated using Ehrhart polynomial\footnote{defined in Section \ref{sec:validity}} in this case, for reference. We observe that while \eqref{CDS_bound} holds for any rate $M^{1-\frac{1}{p}}$, the attached coefficient is very large with small $p$ and strong rates are not effective except at high cardinality; \eqref{BNTT_bound} is slightly above \eqref{CDS_bound}, and both of them show considerable discrepancy with the exact truncation error, verifying Stechkin inequality is not sharp. Estimate \eqref{compare:est1}, on the other hand, is close to the true value, even with $\varepsilon$ large. Besides, the actual minimum cardinality for the estimate to hold is shown {  as}
$M_{\varepsilon} \simeq 1$ for $\varepsilon = 4.0$, 
$M_{\varepsilon} \simeq 10$ for $\varepsilon = 1.0$ and 
$M_{\varepsilon} \simeq 10^3$ for $\varepsilon = 0.3$. 
Also notice that $(N!)^{1/N} \simeq N/e$, \eqref{BNTT_bound} and \eqref{compare:est1} are similar, except for the rate adjusting parameters. While ${1}/({1+\varepsilon})^{1/N}$ in \eqref{compare:est1} can be close to 1, $\xi$ is bounded by $(e-1)/e \simeq 0.65$, resulting in the best convergence attainable by \eqref{BNTT_bound} approximately $ M \exp \left( - 0.65 \left({\frac{M}{|\mathcal{P}|}}\right)^{1/N} \right)$.      

\begin{figure}[h]
\centering
\includegraphics[height=6.0cm]{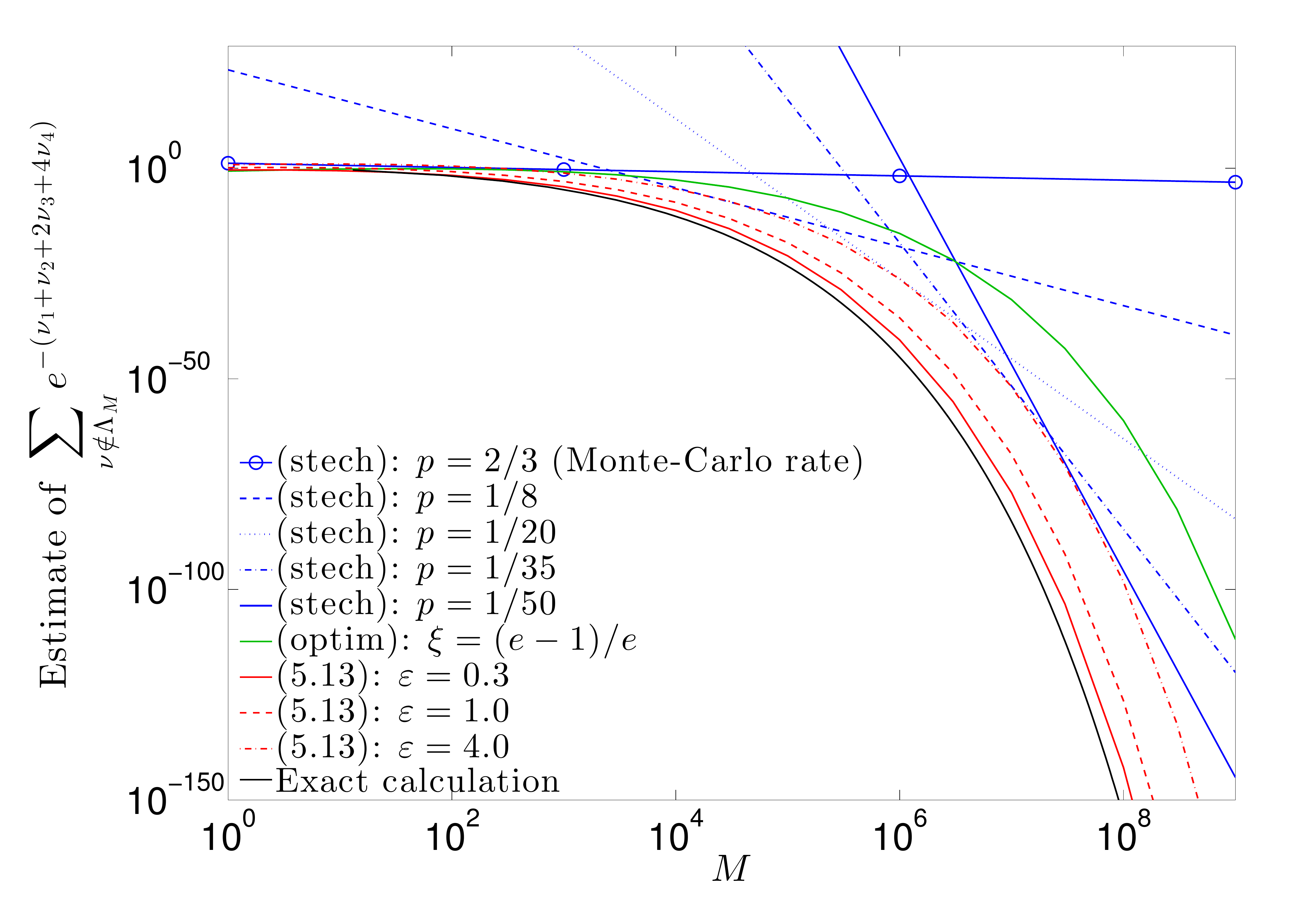} 
\caption{A comparison of our error estimate in computing the series $\sum_{{\bm \nu} \in \mathcal{S}} e^{-(\nu_1 + \nu_2 + 2\nu_3 + 4\nu_4)}$ with those resulting from some previous approaches.}
\label{compareBound}
\end{figure}

We consider next the problem of finding a tight upper bounds of 
\begin{align*}
\textit{error} \, :=\, \sup_{{\bm y}\in\Gamma} \left\| u({\bm y}) - \!\sum_{{\bm \nu} \in \Lambda_M}  t_{\bm \nu}{\bm y}^{\bm \nu} \right\|_{V(D)},
\end{align*} 
assuming $u({\bm z})$ is a holomorphic function in an open neighborhood of the polydisc $\mathcal{O}_{\bm \rho}$ with $\rho_1 = \ldots = \rho_N >1 \ \, \forall N$. We note that \eqref{compare:est1} holds here, since the exponential decay \eqref{compare:eq1} occurs (see Section \ref{estim_coef}), with $\lambda_i = \log\rho_i =: \lambda,\, \forall i$. An isotropic estimate introduced in \cite{BNTT14}, when applied to this error, gives
\begin{align}
\tag{optim-b}
\label{isoOpt}
\textit{error} \le\,  (1-e^{-\lambda/2})^{-N} \exp\left({\frac{\lambda N}{2 e}\log\left(1 - \epsilon \right)} \sqrt[N]{M}\right), 
\end{align}  
where $\epsilon =  \frac{e-1}{e}\left(1 - \frac{1.09}{\sqrt[N]{M}}\right)$. This bound is obtained based on an optimization of a Stechkin-type estimation, also presented in \cite{BNTT14},
\begin{align}
\tag{stech-b}
\label{isoSte}
\textit{error} \le \,   (1-e^{-\lambda/2})^{-N} M^{-1/p} (1- e^{-p\lambda/2})^{-N/p},
\end{align}
for $p > 0$. Another nice result due to \cite{BBL02,BNTT14}, employing complex analysis technique, proves 
\begin{align*}
\textit{error} \le \,  \frac{1}{e^\lambda -1} e^{-\lambda J}, 
\end{align*}
for $M = \left( \begin{array}{c}
\! N + J   \! \\ \! J \! \end{array} \right)$, which implies 
\begin{align}
\label{complex}
\tag{complex}
\textit{error} \le \, \frac{1}{e^\lambda -1}\exp\left(-\lambda (MN!)^{1/N}\right) 
\end{align}
in asymptotic regime. 

Figure \ref{fig:isoBound} plots estimate \eqref{compare:est1} and the upper bounds listed above in case $\lambda = 1$ and $N =8$. 
The exact truncation error in computing the series $\sum_{{\bm \nu} \in \mathcal{S}} \exp(-{\sum_{i=1}^8 \nu_i})$ is also shown. It is interesting to see the \eqref{isoOpt} curve is almost tangent to the \eqref{isoSte} lines, elucidating that \eqref{isoOpt} is obtained by an optimization of \eqref{isoSte}. Again, estimate \eqref{compare:est1} exhibits a much better approximation of the exact truncation error than \eqref{isoSte} and \eqref{isoOpt}. It should, however, be noted that \eqref{isoOpt} is proved to hold with relatively small cardinalities ($M > 1.09^N$), which are not covered by our analysis. The best convergence rate here is given by \eqref{complex}. The advantage of \eqref{complex} lies in the fact that unlike other approaches, it seeks to approximate the remainder of Taylor series $ \left\| u({\bm y}) - \!\sum_{{\bm \nu} \in \Lambda_M}  t_{\bm \nu}{\bm y}^{\bm \nu} \right\|_V$ directly without using triangle inequality. Figure \ref{fig:isoBound} shows a discrepancy between \eqref{complex} and exact calculation of $\sum_{{\bm \nu} \in \mathcal{S}} \exp(-{\sum_{i=1}^8 \nu_i})$, revealing triangle inequality is not sharp in all cases. We are, unfortunately, not aware of an extension of \eqref{complex} outside the isotropic setting.   

\begin{figure}[h]
\centering
\includegraphics[height=6.0cm]{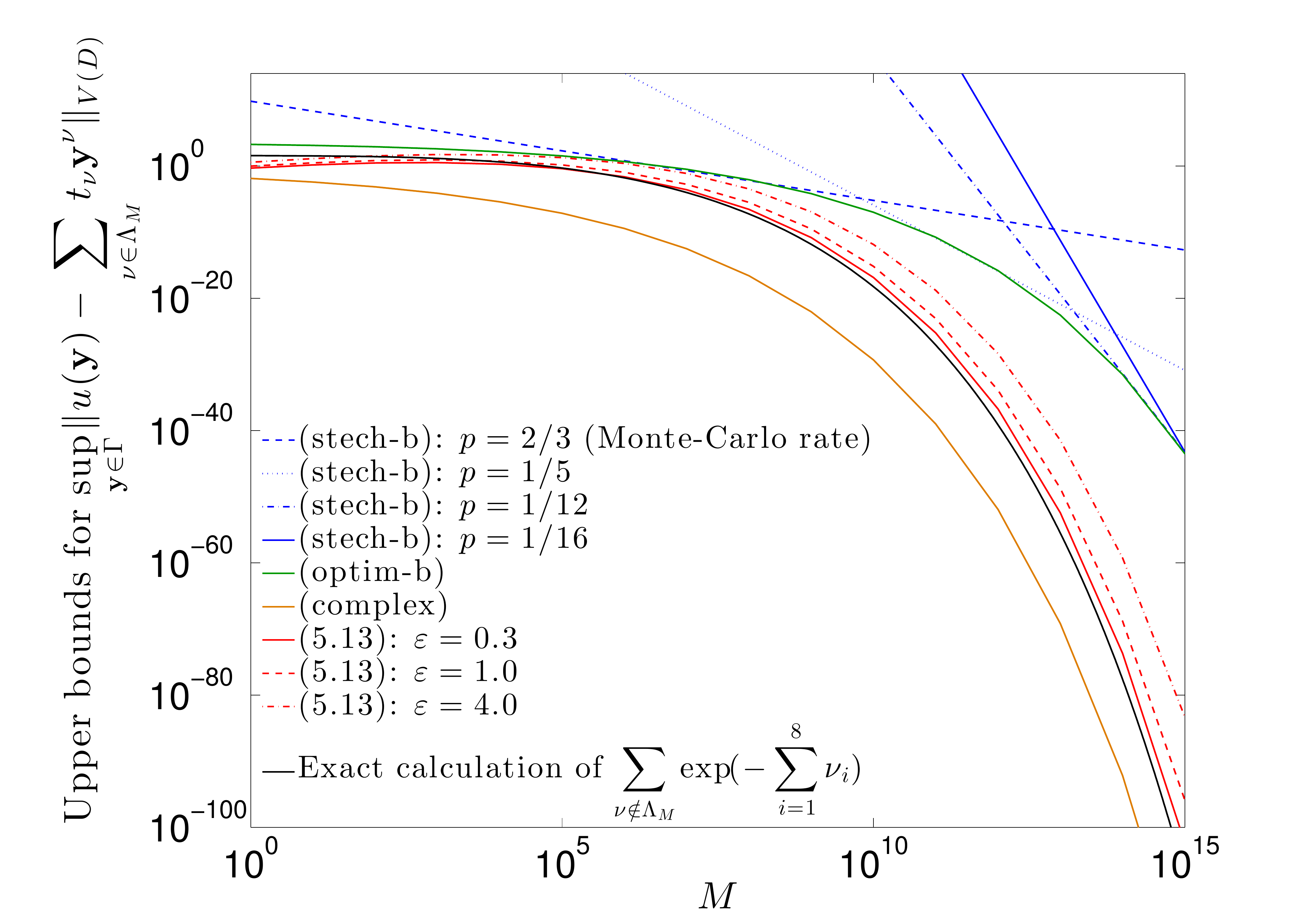} 
\caption{A comparison of our error estimate with those resulting from some previous approaches in an isotropic setting. }
\label{fig:isoBound}
\end{figure}

\section{The optimality of our proposed estimation and pre-asymptotic error analysis: a simplified case}
\label{sec:validity}

In this section, we consider the particular case in which 
\begin{enumerate}[\ \ \ i)]
\label{rat.poly}
\item $\mathcal{P}$ is a rational convex polytope, 
\item ${\mathcal{P}_\tau} = {\tau}  \mathcal{P}$ for all $\tau\in(0,\infty)$.
\end{enumerate} 
This setting, arising from the multi-indexed sequence $(e^{-b({\bm \nu})})_{{\bm \nu} \in \mathcal{S}} $ with 
\begin{align}
\label{b.rat.poly}
b(\bm \nu) =  \sup_{{\bm \lambda} \in {\bm A}}(\sum_{i=1}^N \lambda_i \nu_i ),\mbox{ where } {\bm A}\mbox{ is a \textit{finite} subset of }(\mathbb{Q}^+)^N,
\end{align}
is appropriate for Taylor coefficient estimate of the form \eqref{taylor_est} and to some extend, \eqref{bound_opt} ({  as} discussed in Section \ref{section:conv_anal}). The advantage here is that the number of integer points $\# (\mathcal{P}_j \cap \mathbb{Z}^N)$ can be represented by a computable \textit{Ehrhart quasi-polynomial} of degree $N$ in $j$ (see \cite{BR07}, Chapter 3 and \cite{Sta97}, Chapter 4). In other words, there exist a period $q$ and polynomials $E_0,\ldots,E_{q-1}$ of degree $N$ with leading coefficient $|\mathcal{P}|$ such that $\# (\mathcal{P}_j \cap \mathbb{Z}^N) = E_i(j)$ if $j\equiv i \mod q$. We exploit this property for two tasks: first, to establish a lower bound of $\sum_{{\bm \nu}\notin \Lambda_M} e^{-b({\bm \nu})}$ and verify the sharpness of estimate \eqref{conv_rate1:eq}; second, to calculate the minimum cardinalities $M_{\varepsilon}$ and $M'_{\varepsilon}$ for \eqref{conv_rate1:eq} and \eqref{eq_new:1} to hold (via the computations of the Ehrhart quasi-polynomials) and study the relation between them and the convergence rate. To circumvent the constraints on $M_{\varepsilon}$ and $M'_{\varepsilon}$, an estimate of the truncation errors in the pre-asymptotic regime will be derived. 

\subsection{Lower bound of the truncation errors}
We begin this section with an additional assumption on $b$, which is fulfilled by $b({\bm \nu})$ defined in \eqref{b.rat.poly}. 

\begin{assumption}[Monotonically increasing]
\label{mono_decreasing}
 $b:[0,\infty)^N \to\mathbb{R}$ satisfies: $\forall{\bm \nu},\bm {\mu}\in [0,\infty)^N$, if $\bm{\nu} \le \bm{\mu }$, then $b({\bm \nu}) \le b({\bm \mu})$.   
\end{assumption}
Given this monotone property, the number of integer points inside a superlevel set $\mathcal{P}_\tau$ is always larger than its Lebesgue measure. This observation is verified in the following lemma. 

\begin{lemma}
\label{point_count3}
Assume that $b:[0,\infty)^N \to\mathbb{R}$ is continuous and satisfies Assumption \ref{mono_decreasing}. For $\tau\in (0,\infty)$, denote $\mathcal{P}_\tau = \left\{\bm{\nu} \in [0,\infty)^N:\, b(\bm{\nu})\le \tau \right\}$. We have
\begin{align*}
\# (\mathcal{P}_\tau \cap \mathbb{Z}^N) \ge |\mathcal{P}_{\tau}|,\ \ \forall \tau> 0. 
\end{align*}
\end{lemma}

\sproof
We consider a partition of $[0,\infty)^N$ by the family of cells $(I_{\bm \nu})_{\bm {\nu}\in \mathcal{S}}$ defined as 
\begin{align*}
I_{\bm \nu} = \bigotimes_{1\le i \le N} [\nu_i,\nu_i+1).
\end{align*}  
Denoting $\mathcal{S}^* = \{{\bm \nu} \in \mathcal{S}: \mathcal{P}_{\tau}\cap I_{\bm \nu}\ne \varnothing \}$. If ${\bm \nu} \in \mathcal{S}^*$, by definition, there exists ${\bm \mu} \in I_{\bm \nu}$ such that $b({\bm \mu})\le \tau$. Since ${\bm \nu} \le {\bm \mu}$ and $b$ satisfies Assumption \ref{mono_decreasing}, it gives $b({\bm \nu}) \le b({\bm \mu}) \le \tau$. We have ${\bm \nu} \in \mathcal{P}_{\tau}\cap \mathbb{Z}^N$, which implies $\mathcal{S}^*\subset  \mathcal{P}_{\tau}\cap \mathbb{Z}^N$ and $\#(\mathcal{S}^*)\le \#(\mathcal{P}_{\tau}\cap \mathbb{Z}^N) $. 

On the other hand, there holds 
\begin{align*}
|\mathcal{P}_{\tau}| = \sum_{{\bm \nu}\in\mathcal{S}} |\mathcal{P}_{\tau}\cap I_{\bm \nu}| \le \sum_{{\bm \nu}\in\mathcal{S}^*} |I_{\bm \nu}| = \#(\mathcal{S}^*). 
\end{align*}

We obtain $ |\mathcal{P}_{\tau}| \le \#(\mathcal{S}^*) \le \# (\mathcal{P}_\tau \cap \mathbb{Z}^N) $, as desired. 
\fproof

Now, we proceed to establish a lower bound for the truncation errors of series $\sum_{{\bm \nu}\in \mathcal{S}} e^{-b({\bm \nu})}$ with $b({\bm \nu})$ having the form \eqref{b.rat.poly}. 

\begin{theorem}
\label{prop:lower_bound}
Consider the multi-indexed series $\sum_{{\bm \nu}\in \mathcal{S}} e^{-b({\bm \nu})}$ with $b(\bm \nu)$ given by \eqref{b.rat.poly}. There exists a constant $M^* > 0 $ such that
\begin{align}
\sum_{{\bm \nu}\notin \Lambda_M} e^{-b({\bm \nu})}  \ge  C_{\ell} \, M^{1-\frac{1}{N}} \exp\left(-\left(\frac{M}{|\mathcal{P}| }\right)^{1/N}\right)
\label{lower_bound}
\end{align}
for all $M>M^*$. Here, $\mathcal{P}$ is defined as in \eqref{define.P}, $C_{\ell} = \frac{1}{2}\left(\frac{2}{3}\right)^{1-\frac{1}{N}}  \frac{N|\mathcal{P}|^{\frac{1}{N}} q}{e^q -1} $ where $q$ is the period of Ehrhart quasi-polynomial of $\mathcal{P}$. 
\end{theorem}

\sproof
It is easy to see that $b$ satisfies Assumption \ref{sec4:assump1}. Particularly, $b(\tau{\bm \nu}) = \tau b({\bm \nu})$, $H_{\bm \nu}$ is constant and $\frac{1}{\tau}\mathcal{P}_{\tau} = \mathcal{P} = \left\{\bm{\nu} \in [0,\infty)^N:\, b(\bm{\nu})\le 1 \right\}$ for all $\tau \in (0,\infty),\, {\bm \nu} \in [0,\infty)^N$. By definition of $b$, $\mathcal{P}$ is a rational convex polytope. We can find $q\in \mathbb{N}$ and an $N$-order polynomial $E$ with leading coefficient $|\mathcal{P}|$ such that 
\begin{align}
\label{lower_bound:eq0}
\# (\mathcal{P}_{jq} \cap \mathbb{Z}^N) = E(jq),\ \ \forall j\in  \mathbb{N}.
\end{align}
For $\Lambda_M = \mathcal{P}_{Jq} \cap \mathbb{Z}^N$, it gives 
\begin{gather}
\label{lower_bound:eq1}
\begin{aligned}
& \sum_{{\bm \nu}\notin \Lambda_M} e^{-b({\bm \nu})}  \ge \sum_{j \ge  J} ( \# (\mathcal{P}_{(j+1)q} \cap \mathbb{Z}^N) -  \# (\mathcal{P}_{jq} \cap \mathbb{Z}^N))e^{-(j+1)q} 
\\
 = \, & \sum_{j \ge  J} \left( E(jq+q) - E(jq)\right)e^{-(j+1)q} 
\end{aligned}
\end{gather}
Denoting $E(t) = |\mathcal{P}| t^N + \sum_{i=0}^{N-1} c_i  t^i  ,\ \forall t\in \mathbb{R}$, we have 
\begin{align}
\label{lower_bound:eq2}
E(jq+q) - E(jq) \ge q |\mathcal{P}| N (jq)^{N-1} -  q  \sum_{i=0}^{N-1} |c_i| i  (jq + q)^{i-1}.
\end{align}
There exists $\Upsilon_1 > 0$ satisfying
\begin{align}
\label{lower_bound:eq3}
  \sum_{i=0}^{N-1} |c_i| i  (jq + q)^{i-1} \le \frac{1}{2}  |\mathcal{P}| N (jq)^{N-1},\ \  \forall j\in \mathbb{N},\, j > \Upsilon_1. 
\end{align}
Combining \eqref{lower_bound:eq1}--\eqref{lower_bound:eq3} yields for $J > \Upsilon_1$,
\begin{gather}
\label{lower_bound:eq4}
\begin{aligned}
& \sum_{{\bm \nu}\notin \Lambda_M} e^{-b({\bm \nu})}  \ge \frac{1}{2} q|\mathcal{P}| N \sum_{j \ge  J} (jq)^{N-1} e^{-(j+1)q}
\\
\ge\, &  \frac{1}{2} N q^N J^{N-1}  |\mathcal{P}|  \sum_{j \ge  J} e^{-(j+1)q} = \frac{1}{2} Nq (qJ)^{N-1}  |\mathcal{P}| \frac{e^{-qJ}}{e^q-1}. 
\end{aligned}
\end{gather} 
We need to write this estimate in term of the cardinality $M$. First, notice that $b$ satisfies Assumption \ref{mono_decreasing}, there holds
\begin{align}
\label{lower_bound:eq4b}
|\mathcal{P}| (Jq)^{N} = |\mathcal{P}_{Jq}|  \le  \# (\mathcal{P}_{Jq} \cap \mathbb{Z}^N).
\end{align} 
Applying Theorem \ref{conv_rate1}, it gives $ |\mathcal{P}| = \lim\limits_{j\to\infty} \frac{1}{(jq)^{N}}\cdot \# (\mathcal{P}_{jq} \cap \mathbb{Z}^N) $. We can choose $ \Upsilon_2>0$ such that for all $j\in \mathbb{N},\, j > \Upsilon_2$, 
\begin{gather}
\begin{aligned}
\label{lower_bound:eq5}
-\frac{1}{2}|\mathcal{P}| \le \, |\mathcal{P}| - \frac{1}{(jq)^{N}}\cdot \# (\mathcal{P}_{jq} \cap \mathbb{Z}^N).
\end{aligned}
\end{gather}
Since $M =  \# (\mathcal{P}_{Jq} \cap \mathbb{Z}^N)$, from \eqref{lower_bound:eq4b} and \eqref{lower_bound:eq5}, one has 
\begin{align}
|\mathcal{P}| (Jq)^{N} \le M  \le \frac{3}{2} |\mathcal{P}|(Jq)^{N}\ \mbox{ for $J>\Upsilon_2$}.
\label{lower_bound:eq6}
\end{align}
Combining \eqref{lower_bound:eq4} and \eqref{lower_bound:eq6} gives 
\begin{align*}
 \sum_{{\bm \nu}\notin \Lambda_M} e^{-b({\bm \nu})} \ge C_{\ell}   M^{1-\frac{1}{N}} \exp\left(-\left(\frac{M}{|\mathcal{P}| }\right)^{1/N}\right), 
\end{align*}
where $C_{\ell} = \frac{1}{2}\left(\frac{2}{3}\right)^{1-\frac{1}{N}}  \frac{N|\mathcal{P}|^{\frac{1}{N}} q}{e^q -1} $. The proof is now complete. 
\fproof

Theorem \ref{conv_rate1} and Theorem \ref{prop:lower_bound} reveal that for $b({\bm \nu})$ given by \eqref{b.rat.poly}, the asymptotic truncation error of $\sum_{{\bm \nu}\in \mathcal{S}} e^{-b({\bm \nu})}$ can be bounded from below and above as 
\begin{align*}
 C_{\ell} \, M^{1-\frac{1}{N}} \exp\left(-\left(\frac{M}{|\mathcal{P}| }\right)^{1/N}\right) &\le \sum_{{\bm \nu}\notin \Lambda_M} e^{-b({\bm \nu})}\\  
 &\le  C_{u}(\varepsilon)\,  M \exp \left( - \left({\frac{M}{|\mathcal{P}|(1+\varepsilon)}}\right)^{1/N} \right), 
\end{align*}
where $C_{\ell}$ and $C_u(\varepsilon)$ are mild constants in comparison with the total bounds. The optimality of our estimation is verified in these cases.  

\subsection{Asymptotic minimum cardinalities and their relation with the convergence rate}

In this section, we will apply Ehrhart (quasi-)polynomial to investigate the minimum cardinality for our asymptotic convergence rate to hold. Recall that for any $\varepsilon > 0$, the upper estimates \eqref{conv_rate1:eq} and \eqref{eq_new:1} occur with $J> J_{\varepsilon } = \max\left\{\frac{2}{e^{1/N} -1}, \Delta_{\varepsilon} \right\}$ and $J> {J}'_{\varepsilon} =  \max\left\{\frac{1}{e^{1/N} -1}, \Delta_{\varepsilon} \right\}$, respectively. The first constraints in both conditions are straightforward and we focus on quantifying $\Delta_\varepsilon$. From \eqref{eq:3b}, $\Delta_{\varepsilon}$ is the positive real number such that 
\begin{align}
\label{min_card:eq1}
\frac{1}{2} j^N |\mathcal{P}| \le \# (\mathcal{P}_j \cap \mathbb{Z}^N) \le j^N |\mathcal{P}| (1+\varepsilon),\ \forall j\in \mathbb{N},\, j>\Delta_{\varepsilon}.
\end{align} 
In case $\mathcal{P}$ is a rational convex polytope and $\mathcal{P}_{\tau} = \tau\mathcal{P}\ \forall \tau\in (0,\infty)$, we can ignore the left inequality of \eqref{min_card:eq1}, which by Lemma \ref{point_count3} is true for all $j\in \mathbb{N}$. There exists a (quasi-) polynomial 
\begin{align}
\label{min_card:eq2}
E^*(j) = |\mathcal{P}| j^N + \sum_{i=0}^{N-1} c^*_i(j) j^i, 
\end{align}
with $c^*_i:\mathbb{N} \to \mathbb{Q}$ being a periodic function with integer period $q$ such that 
\begin{align}
\label{min_card:eq3}
 \# (\mathcal{P}_j \cap \mathbb{Z}^N)  = E^*(j),\ \forall j\in \mathbb{N},
\end{align}  
see \cite{BR07}, Chapter 3 and \cite{Sta97}, Chapter 4. Replacing \eqref{min_card:eq3} to \eqref{min_card:eq1}, $\Delta_\varepsilon$ can be characterized as the largest among the solutions of 
\begin{align*}
 \varepsilon |\mathcal{P}| j^N - \sum_{i=0}^{N-1} c^*_i(j)  j^i  &= 0.
\end{align*}
The numerical computation of formula of Ehrhart polynomial $E^*(j)$ can be done efficiently \cite{DHTY03}, allowing us to quantify $\Delta_\varepsilon$ and the theoretical minimum cardinality $M_{\varepsilon}$ $(=E^*(J_{\varepsilon}))$ accurately. We present a brief study on the relation between $M_{\varepsilon}$ and $\varepsilon$ for some polytopes, including:
\begin{itemize}
\item (P.1): $b({\bm \nu}) = \sum_{i=1}^4 \nu_i$ ($N=4$),
\item (P.2): $b({\bm \nu}) = \nu_1 + \nu_2 + 2\nu_3 + 4\nu_4$ ($N=4$),
\item (P.3):  $b({\bm \nu}) = \sum_{i=1}^8 \nu_i$ ($N=8$),
\item (P.4): $b({\bm \nu}) = \sum_{i=1}^8  \frac{\nu_i}{2^{i-3}} $ ($N=8$),
\item (P.5): $b({\bm \nu}) = \sup\left\{\frac{1}{2}\sum_{i=1}^8  {\nu_i},\, \frac{5}{16}\sum_{i=1}^8  {\nu_i} + \frac{5}{16}  {\nu_j} : 1\le j \le 8  \right\} $ ($N=8$),
\item (P.6): $b({\bm \nu}) = \sup\left\{\frac{1}{5}\sum_{i=1}^8  {\nu_i},\, \frac{1}{8}\sum_{i=1}^8  {\nu_i} + \frac{1}{8}  {\nu_j} : 1\le j \le 8  \right\} $ ($N=8$). 
\end{itemize} 
(P.1)-(P.4) correspond to $4$- and $8$-simplices with different levels of anisotropy. The lengths of edges connecting the origin and other vertices are equal for (P.1) and (P.3) and slightly vary for (P.2), while (P.4) is quite a skinny simplex. On the other hand, (P.5) is a truncated, enlarged version of (P.3) where the vertices are at $\frac{1}{5}$ of the way along the axis edges and $\frac{2}{5}$ along other edges, resulting in a polytope with $65$ vertices. (P.6) in turn is obtained through an enlargement of (P.5). {  We note that (P.2) and (P.3) correspond to the coefficient bounds illustrated in Section \ref{sec:comparison}}.

Figure \ref{Mmin} shows the variation of $M_{\varepsilon}$ and $M'_{\varepsilon}$ as well as the \textit{rate adjusting parameter} $1/(1+\varepsilon)^{1/N}$ in the estimates \eqref{conv_rate1:eq} and \eqref{eq_new:1} with respect to $\varepsilon$. The other parameter $C_u(\varepsilon)$ is negligible except for $\varepsilon$ very large and not plotted here. The formulas of Ehrhart polynomials are calculated using the software package \texttt{LattE} \cite{BBDD+13}. First, we observe that choosing a smaller $\varepsilon$ gives a stronger convergence, yet $M_{\varepsilon}$ must also be increased. The good news is that while the best convergence $M \exp \left( - \left({\frac{M}{|\mathcal{P}|}}\right)^{1/N} \right)$ is realized only as $\varepsilon \to 0$, $\varepsilon$ need not to be small to obtain a strong rate, especially in high dimension. For instance, $\varepsilon = 1.0$ gives the rate $\sim M \exp \left( -{0.83}\left(\frac{M}{|\mathcal{P}|}\right)^{1/N} \right)$ with $N = 4$ and $\sim M \exp \left( - 0.92 \left({\frac{M}{|\mathcal{P}|}}\right)^{1/N} \right)$ with $N = 8$. 
 
Not surprisingly, $M_{\varepsilon}$ and $M'_{\varepsilon}$ is shown to be larger for higher dimension. For a fixed $N$, the anisotropy of the polytopes significantly impacts $M_{\varepsilon}$ and $M'_{\varepsilon}$: these values are close for (P.3) and (P.5), which possess different shapes and scales but span equally in coordinate axes, and much larger for (P.4), the simplex with skinny shape. Generally, increasing $\varepsilon$ alleviates the restrictions on $M_{\varepsilon}$ and $M'_{\varepsilon}$, as this will reduce $\Delta_{\varepsilon}$. The strategy is, however, ineffective once $\frac{2}{e^{1/N}-1}$ (or $\frac{1}{e^{1/N}-1}$) exceeds $\Delta_\varepsilon$ and dominates \eqref{J_cond} and \eqref{J'_cond}, at which point, these conditions can no further be relaxed. Thus, while $M_{\varepsilon}$ and $M'_{\varepsilon}$ are almost not affected by the scale of polytopes with $\varepsilon$ close to $0$, their lower bounds (imposed by $J\gtrsim \frac{1}{e^{1/N} -1} \simeq N$) are more restrictive for large polytopes; in such cases, mild constraints on $M_{\varepsilon}$ and $M'_{\varepsilon}$ may be unattainable. This fact is illustrated by a comparison of two similar polytopes (P.5) and (P.6) in Figure \ref{Mmin}: $M_{\varepsilon}$ and $M'_{\varepsilon}$ eventually stop to decay in both cases, but the bound is higher for (P.6), the polytope with larger scale.   

In short, our asymptotic convergence analysis applies to the range $J\ge N$. In the next part, we propose an alternative estimate of truncation errors, which is effective in the pre-asymptotic regime $J<N$. {  Let us remark that the actual condition on $M_{\varepsilon}$ for \eqref{conv_rate1:eq} to hold can be much milder than the theoretical minimum cardinality posed by Theorem \ref{conv_rate1} and investigated here, as shown in Figure \ref{compareBound} and \ref{fig:isoBound}.}

\begin{figure}[h]
\centering
\includegraphics[height=3.5cm]{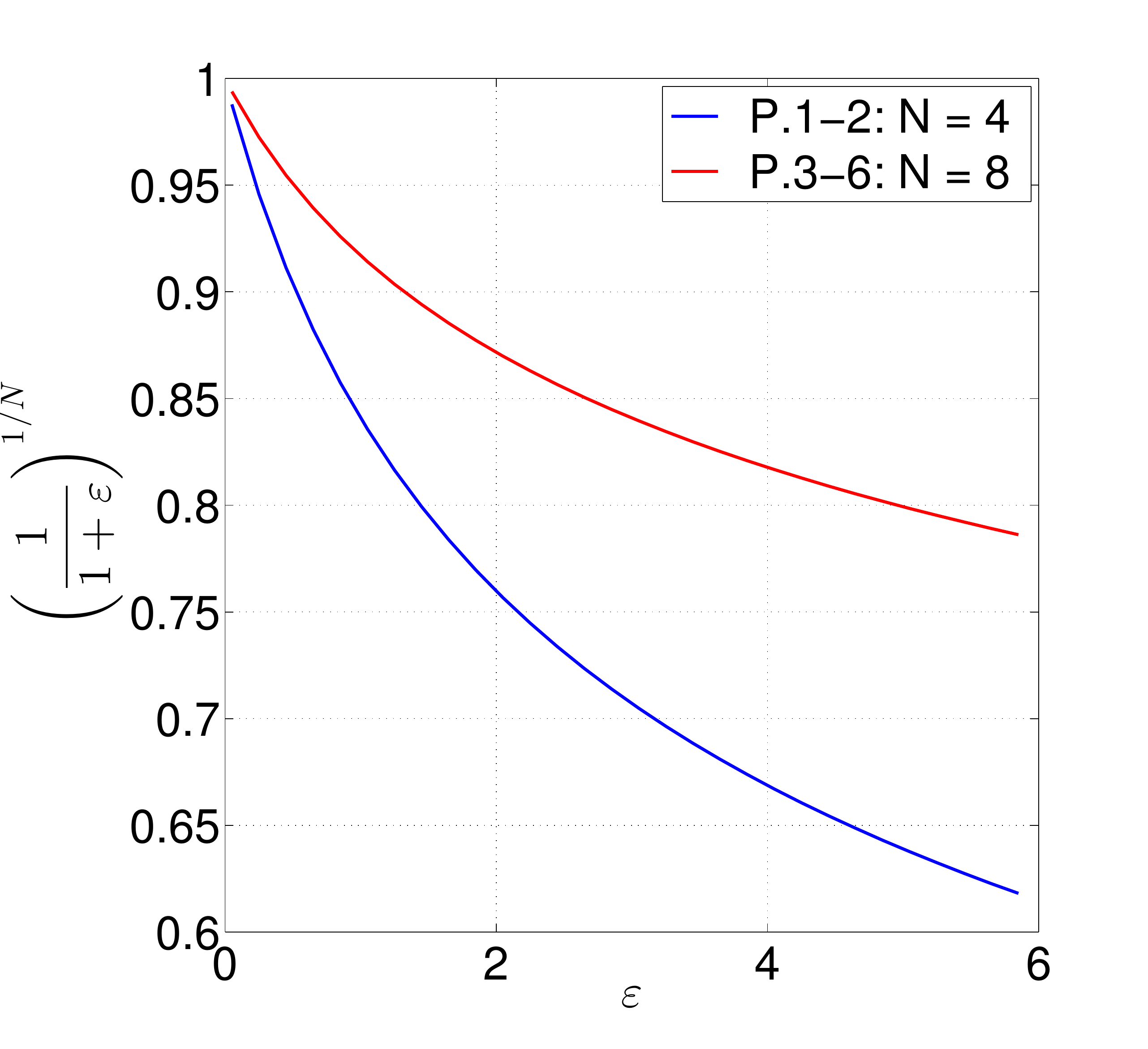} 
\includegraphics[height=3.5cm]{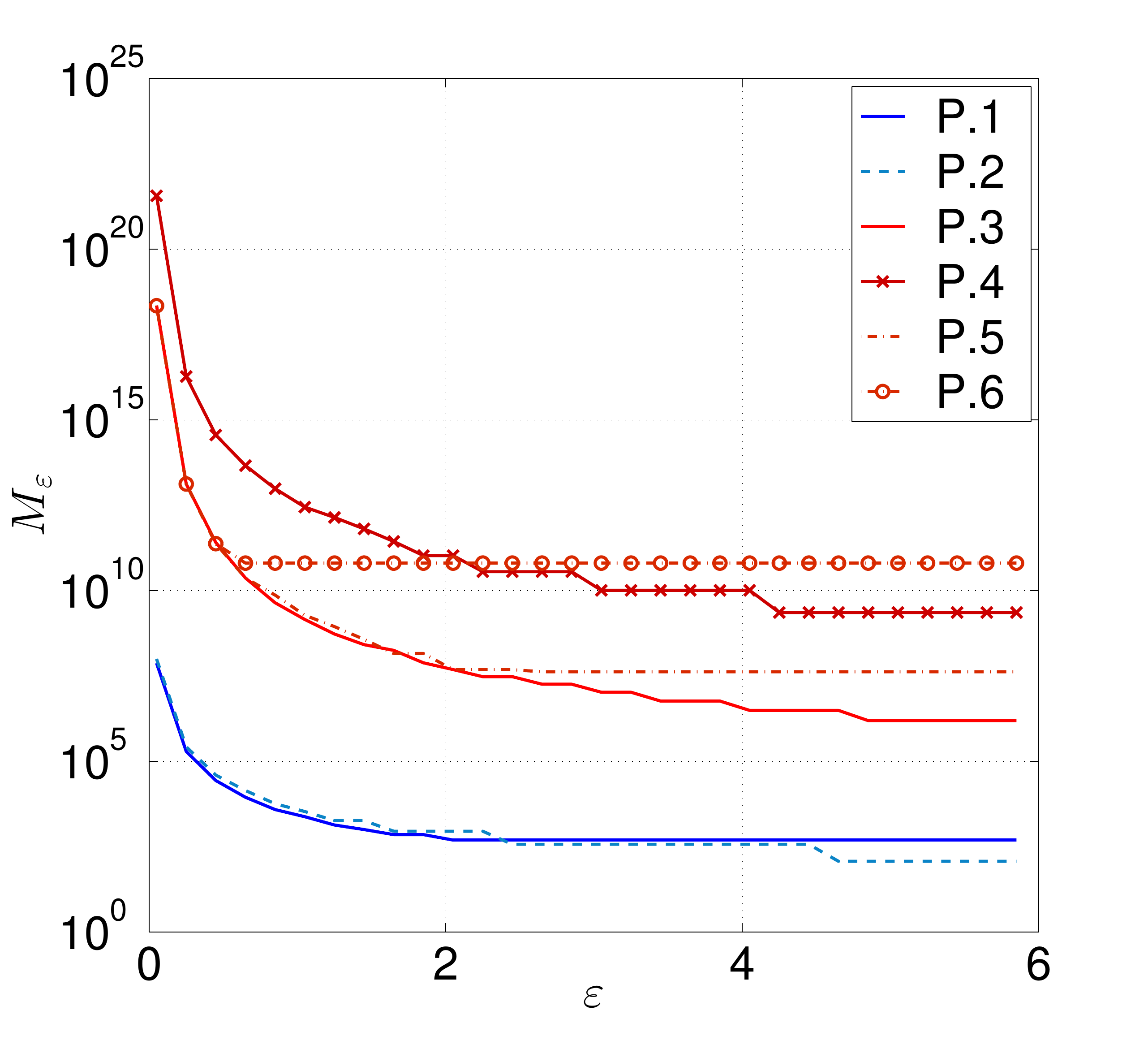} 
\includegraphics[height=3.5cm]{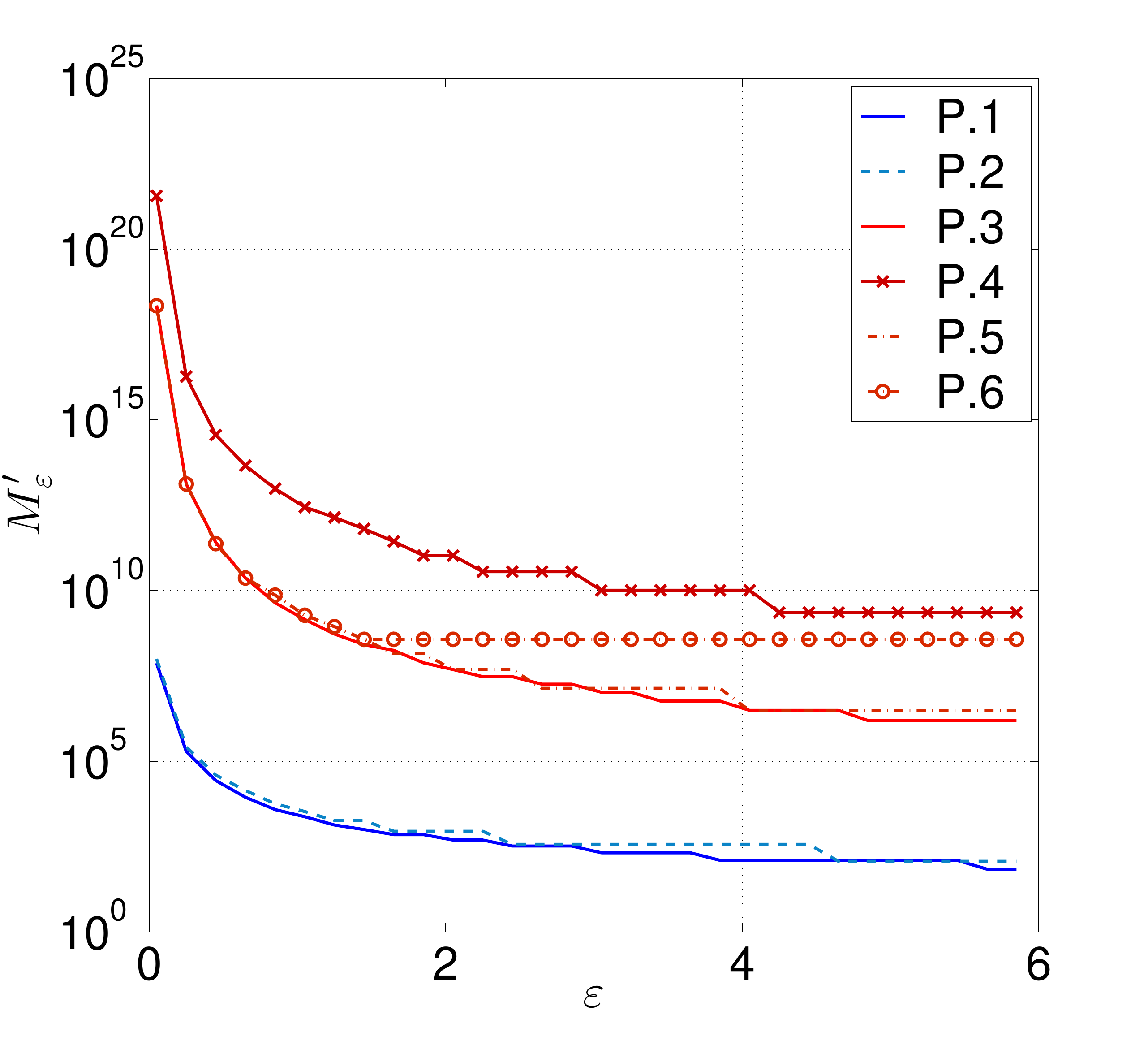} 
\caption{The variation of the rate adjusting parameter and theoretical minimum cardinalities $M_{\varepsilon}$ and $M'_{\varepsilon}$ for the upper estimate \eqref{conv_rate1:eq} with respect to $\varepsilon$.}
\label{Mmin}
\end{figure}

\subsection{A pre-asymptotic estimate of truncation errors} 
\label{preasymp}
To acquire an estimation of $\sum_{{\bm \nu} \notin \Lambda_M} e^{-b({\bm \nu})}$ in pre-asymptotic regime, following the arguments in Theorem \ref{conv_rate1}, non-asymptotic bounds of $\#(\mathcal{P}_j \cap \mathbb{Z}^N)$ and $\sum_{j\ge J} j^N e^{-j}$ need to be established. 
An upper bound of $\#(\mathcal{P}_j \cap \mathbb{Z}^N)$ is derived in the following lemma. 
\begin{lemma}
\label{point_count4}
Let $b:[0,\infty)^N \to \mathbb{R}$ be continuous and satisfy Assumption \ref{mono_decreasing}. Assuming that $b(\tau{\bm \nu}) = \tau b({\bm \nu})$ for all $\tau \in (0,\infty), {\bm \nu} \in [0,\infty)^N$. For $\tau \in (0,\infty)$, denote $\mathcal{P}_{\tau} = \{{\bm \nu}\in [0,\infty)^N: b({\bm \nu}) \le \tau\}$. There follows 
\begin{align*}
\#(\mathcal{P}_{j}\cap \mathbb{Z}^N) \le j^N\cdot  \#\left(\mathcal{P} \cap \mathbb{Z}^N\right),\ \ \forall j \in \mathbb{N},
\end{align*}
where $\mathcal{P} =\{{\bm \nu}\in [0,\infty)^N: b({\bm \nu}) \le 1\}$ ( $= \frac{1}{\tau}\mathcal{P}_\tau$ for all $\tau$). 
\end{lemma}

\sproof
Since $b(\tau{\bm \nu}) \equiv \tau b({\bm \nu})$, we have $\mathcal{P}_\tau = \tau \mathcal{P},\, \forall \tau >0$, thus, 
\begin{align*}
\#(\mathcal{P}_{j}\cap \mathbb{Z}^N) = \#(j \mathcal{P}\cap \mathbb{Z}^N) =  \#\left(\mathcal{P}\cap \frac{1}{j}\mathbb{Z}^N\right), \ \forall j\in \mathbb{N}.
\end{align*}

Given ${\bm \mu} \in \left(\mathcal{P}\cap \frac{1}{j}\mathbb{Z}^N\right) $, $\bm \mu$ can be written uniquely in the form 
\begin{align*}
{\bm \mu} = {\bm \nu} + \bigotimes_{i=1}^N \frac{r_i}{j},
\end{align*}
where $\bm \nu\in \mathcal{S}$ and $r_i\in \mathbb{Z}, 0\le r_i \le j-1,\, \forall\, 1\le i \le N$. 

Since ${\bm \nu} \le {\bm \mu}$ and $b$ satisfies Assumption \ref{mono_decreasing}, it gives $b({\bm \nu}) \le b({\bm \mu}) \le 1$ and, consequently, ${\bm \nu}\in \mathcal{P} \cap \mathbb{Z}^N$. We have 
\begin{align*} 
 &\#\left(\mathcal{P}\cap \frac{1}{j}\mathbb{Z}^N\right) \\
 \le & \#  \left\{{\bm \nu} + \bigotimes_{i=1}^N \frac{r_i}{j}:{\bm \nu}\in \mathcal{P} \cap \mathbb{Z}^N , r_i\in \mathbb{Z}, 0\le r_i \le j-1,\, \forall\, 1\le i \le N\right\} 
 \\
  = & j^N\cdot  \#\left(\mathcal{P} \cap \mathbb{Z}^N\right),
\end{align*}
as desired. 
\fproof

Next, we give a non-asymptotic estimate of $\sum_{j\ge J} j^N e^{-j}$ based on tight approximation of $\sum_{j\le J-1} j^N e^{-j}$ for $J\le N+1$. Indeed, since $\tau\mapsto  \tau^N e^{-\tau}$ is increasing in $[0,N]$, we have
\begin{align}
\label{pre_asym:eq0}
 \sum_{j =1}^{J-1} j^N e^{-j} \ge  \int_{0}^{J-1} \tau^N e^{-\tau}d\tau.
\end{align}
Applying Theorem 4.1, \cite{Neu13}, yields
\begin{align}
\label{pre_asym:eq0b}
\int_{0}^{J-1} \tau^N e^{-\tau}d\tau \ge  \frac{(J-1)^{N+1}}{N+1} \exp\left(-\frac{(J-1)(N+1)}{N+2}\right).
\end{align}
Combining \eqref{pre_asym:eq0} and \eqref{pre_asym:eq0b}, it gives
\begin{align}
\sum_{j\ge J} j^N e^{-j} &= \sum_{j =1}^{\infty} j^N e^{-j} -  \sum_{j =1}^{J-1} j^N e^{-j} \notag
\\
& \le \sum_{j =1}^{\infty} j^N e^{-j}  -  \frac{(J-1)^{N+1}}{N+1} \exp\left(-\frac{(J-1)(N+1)}{N+2}\right). \label{pre_asym:eq0c}
\end{align}
A mathematical formula of the sum $ \sum_{j =1}^{\infty} j^N e^{-j}  $ is not accessible. However, it is independent of $J$ and can be written in term of the well-studied \textit{polylogarithm functions}
\begin{align}
\label{polylogarithm}
Li_s(z) = \sum_{j=1}^\infty \frac{z^j}{j^s},\ \mbox{ for }z\in \mathbb{C}, |z| < 1, s\in\mathbb{R}.
\end{align} 
see \cite{Lew81,Lew91}. Combining \eqref{pre_asym:eq0c} and \eqref{polylogarithm}, we have proved the following Lemma:
\begin{lemma}
\label{pre_asym:lemma}
For any $N,J \in \mathbb{N}$, if $J\le N + 1$, it gives
\begin{align}
\label{pre_asym:eq1}
\sum_{j\ge J} j^N e^{-j} \le Li_{-N}\left(1/e\right) - \frac{(J-1)^{N+1}}{N+1} \exp\left(-\frac{(J-1)(N+1)}{N+2}\right).
\end{align}
\end{lemma}

In Figure \ref{pre_asymp:fig}, we compare the performance of the asymptotic bound \eqref{lemma:eq3} and the pre-asymptotic bound \eqref{pre_asym:eq1} in estimating the truncation error of $\sum_{j = 1}^{\infty} j^N e^{-j}$ for $N = 20$. The pre-asymptotic estimate shows an excellent agreement with true value for small $J$; however, it cannot capture the error decay when $J$ is big. The asymptotic bound, on the other hand, successfully predicts the convergence rate of $\sum_{j = 1}^{\infty} j^N e^{-j}$, but is not effective with small $J$. %
\begin{figure}[h]
\centering
\includegraphics[height=5.0cm]{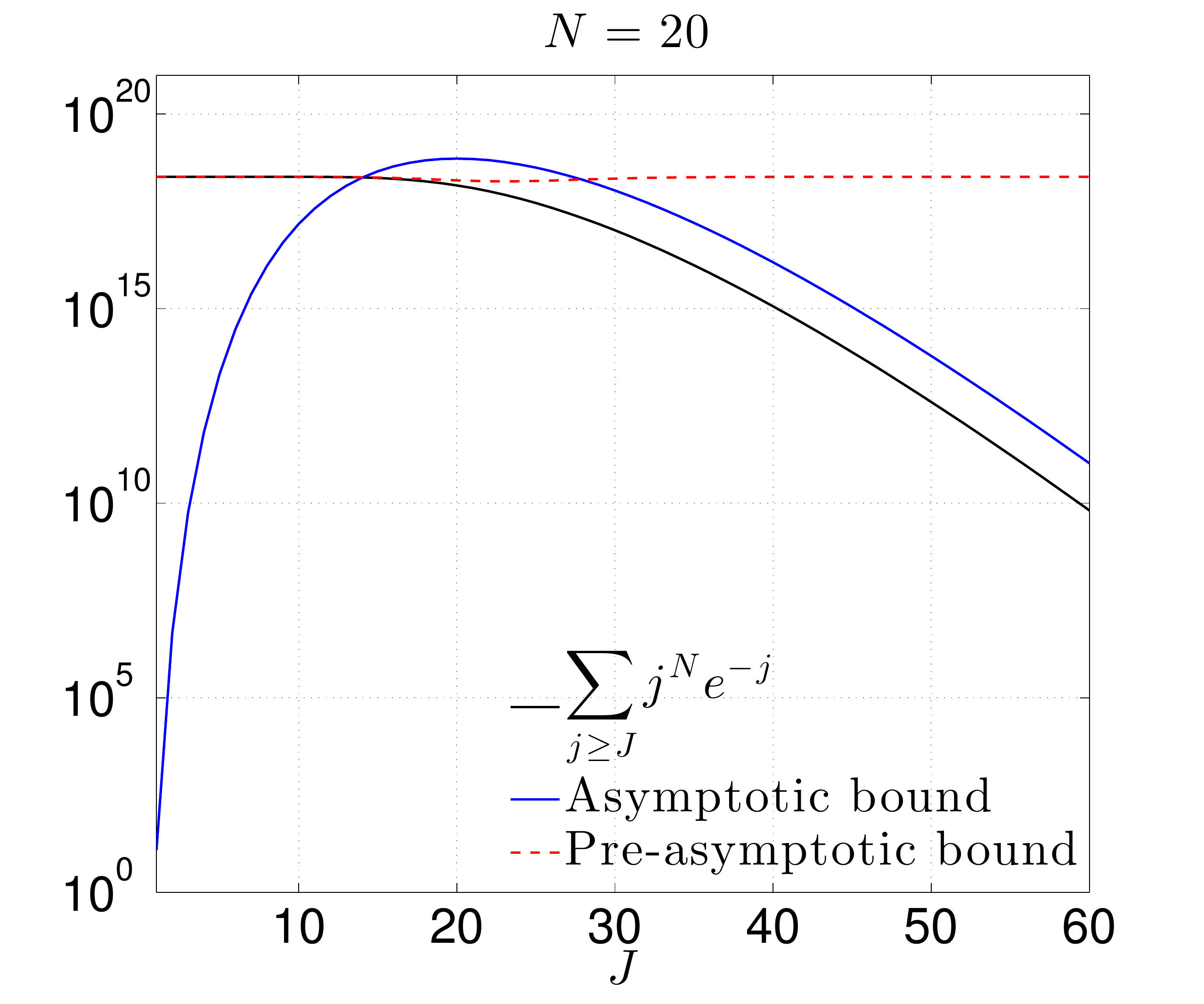} 
\caption{A comparison of the asymptotic bound \eqref{lemma:eq3} and the pre-asymptotic bound \eqref{pre_asym:eq1} in estimating $\sum_{j\ge J} j^N e^{-j}$ for $N = 20$.}
\label{pre_asymp:fig}
\end{figure}

We are now in the position to prove a pre-asymptotic estimation of \\
$\sum_{{\bm \nu}\notin \Lambda_M} e^{-b({\bm \nu})}$. 

\begin{theorem}
Consider the multi-indexed series $\sum_{{\bm \nu}\in \mathcal{S}} e^{-b({\bm \nu})}$ with $b(\bm \nu)$ being continuous and satisfying Assumption \ref{mono_decreasing}. Assuming that $b(\tau{\bm \nu}) = \tau b({\bm \nu})$ for all $\tau \in (0,\infty), {\bm \nu} \in [0,\infty)^N$. For $\tau\in (0,\infty)$, denote $\mathcal{P}_\tau = \left\{\bm{\nu} \in [0,\infty)^N:\, b(\bm{\nu})\le \tau \right\}$ and $\Lambda_M$ the set of indices corresponding to $M$ largest $e^{-b(\bm{\nu})}$. Define $\mathcal{P}$ as in \eqref{define.P}. 

For $M\in \mathbb{N}$, if $M\le \#(\mathcal{P}_N \cap \mathbb{Z}^N)$, there holds 
\begin{align}
\sum_{{\bm \nu}\notin\Lambda_M} e^{-b({\bm \nu})} \le e\sigma \left [Li_{-N}\left(1/e\right) - \frac{1}{N+1}\left(\frac{M}{\sigma} \right)^{\frac{N+1}{N}} \exp\left(-\frac{M^{\frac{1}{N}}(N+1)}{{\sigma}^{\frac{1}{N}}(N+2)}\right) \right ]. 
\label{pre_asym:eq3}
\end{align}

Here, $\sigma = {\#(\mathcal{P} \cap \mathbb{Z}^N)} $ and $Li$ denotes the polylogarithm function. 
\end{theorem}

\sproof
Applying Lemma \ref{point_count4}, it gives 
\begin{align}
\label{pre_asym:eq2}
\#(\mathcal{P}_{j}\cap \mathbb{Z}^N) \le j^N\cdot  \#\left(\mathcal{P} \cap \mathbb{Z}^N\right),\, \forall\, j\in \mathbb{N}.
\end{align}

To estimate $\sum_{{\bm \nu}\notin\Lambda_M} e^{-b({\bm \nu})}$, it is sufficient to consider this sum with $\Lambda_M = \mathcal{P}_J \cap \mathbb{Z}^N,\, J\in \mathbb{N},\, J\le N$. We have 
{\allowdisplaybreaks
\begin{align*}
\sum_{{\bm \nu}\notin \mathcal{P}_J \cap \mathbb{Z}^N }e^{-b({\bm \nu})} \le\, & \sum_{j \ge  J} ( \# (\mathcal{P}_{j+1} \cap \mathbb{Z}^N) -  \# (\mathcal{P}_{j} \cap \mathbb{Z}^N))e^{-j} 
\\
\le\, & \#(\mathcal{P} \cap \mathbb{Z}^N) \sum_{j \ge  J}  {(j+1)}^{N}  e^{-j}
\\
\le\, & e \cdot \#(\mathcal{P} \cap \mathbb{Z}^N) \left [Li_{-N}\left(1/e\right) - \frac{J^{N+1}}{N+1} \exp\left(-\frac{J(N+1)}{N+2}\right) \right ],
\end{align*} 
by applying Lemma \ref{pre_asym:lemma}. 

From \eqref{pre_asym:eq2}, it gives $J\ge \left(\frac{M}{\#(\mathcal{P} \cap \mathbb{Z}^N)} \right)^{1/N}$. It is easy to see that the mapping $j \mapsto \frac{j^{N+1}}{N+1} \exp\left(-\frac{j(N+1)}{N+2}\right)$ is increasing in $[0,N]$. There follows 
\begin{align*}
\sum_{{\bm \nu}\notin\Lambda_M} e^{-b({\bm \nu})} \le e\sigma \left [Li_{-N}\left(1/e\right) - \frac{1}{N+1}\left(\frac{M}{\sigma} \right)^{\frac{N+1}{N}} \exp\left(-\frac{M^{\frac{1}{N}}(N+1)}{{\sigma}^{\frac{1}{N}}(N+2)}\right) \right ],
\end{align*}
}
where $\sigma = {\#(\mathcal{P} \cap \mathbb{Z}^N)} $, implying the assertion \eqref{pre_asym:eq3}. 
\fproof

\begin{remark}
The pre-asymptotic analysis presented above does not employ Ehrhart polynomials, hence applies to a wider class of $b$ than those given by \eqref{b.rat.poly}. In this subsection, $b$ only needs to be continuous, satisfy Assumption \ref{mono_decreasing} and $b(\tau{\bm \nu}) = \tau b({\bm \nu})$ for all $\tau \in (0,\infty), {\bm \nu} \in [0,\infty)^N$.
\end{remark}

\section{Concluding remarks}
\label{sec:conclusion}

We present a new approach for analyzing the convergence of quasi-optimal Taylor and Legendre approximations for parameterized PDEs with {  finite-dimensional} deterministic and stochastic coefficients. {  The analysis also gives an accessible estimates for best $M$-term approximation errors without employing Stechkin inequality.} The advantage of our framework, which is demonstrated through several theoretical examples herein, includes its applicability to a general class of { coefficient decay} and the sharp estimates of asymptotic errors. This work is restricted to linear elliptic equations with input coefficients depending affinely on the parameter. We expect similar results to hold in different settings with finite parametric dimension, particularly nonlinear elliptic PDEs, initial value problems and parabolic equations \cite{CCS14,HS13,HS13b,HS13c}, as our analysis only depends on the polynomial coefficient estimates. 

Developing algorithms for identifying quasi-optimal subspaces is the next natural and essential step. Two potential types of procedures for building the subspaces corresponding to sharp estimates of the coefficients $c_{\bm \nu}$ includes a priori and a posteriori approaches. In the first approach, the estimates for $c_{\bm \nu}$ are derived a priori using knowledge on the input coefficient $a(x,{\bm y})$. Analytical studies reveal that if the complex continuation of $a(x,{\bm y})$ is an analytic function in $\mathbb{C}^N$ then a theoretical decaying rate ${\bm \rho}^{-{\bm \nu}}$ of $c_{\bm \nu}$ (with ${\bm \rho} = (\rho_i)_{1\le i\le N}$ representing the size of certain $N$-dimensional complex domains where real part of $a(x,{\bm y})$ is bounded away from 0) can be proved. The exploration of polynomial subspaces thus reduces to the specification of such domains (or ${\bm \rho}$ in particular), which is expectedly significantly less computational demanding. Recent study \cite{BNTT14} for a priori constructed Total Degree subspace found that while the theoretical estimates were not sharp, they could still provide good prediction on the anisotropy of the index sets. However, in practice, most analytical coefficient bounds lead to subspaces much more complicated than Total Degree and the determination of ${\bm \rho}$ in several cases is nontrivial, see \cite{BNTT14,CCS14,CDS10,CDS11}. It is important to develop, implement and test of effectiveness of a priori algorithms in such settings. 

Research on {\em a posteriori} procedures may be pursued in three directions. The first strategy finds the quasi-optimal index set using the theoretical coefficient estimates, but with ${\bm \rho}$ determined sharply in an a posteriori manner (by exact calculation of the decaying rate of $c_{\bm \nu}$ in each direction $i$, $i=1,\ldots,N$), instead of a priori (by the definition of $a(x,{\bm y})$ as in above). The second strategy adaptively builds nested sequence $(\Lambda_M)_{M\ge 0}$ of quasi-optimal index sets $\Lambda_M$ at a cost that scales linearly in $\#(\Lambda_M)$. Given $\Lambda_M$, we construct $\Lambda_{M+1}$ by enriching $\Lambda_M$ with the most effective indices ${\bm \nu}$ in its neighborhood (denoted by $\mathcal{M}(\Lambda_M)$), which results in the best residual reduction. 
The third strategy first evaluates $u({\bm y})$ on certain finite subset of $\Gamma$ and then constructs the quasi-optimal subspace based on estimates of coefficients $c_{\bm \nu} = \int_\Gamma u({\bm y}) {\bm \Psi}_{\bm \nu}({\bm y}) d{\bm y}$ using non-intrusive methods, e.g., Monte-Carlo, collocation. We expect the exploration cost for this approach, mostly coming from the evaluation of $u({\bm y})$, to be a fraction of cost for computing the solution. 

Finally, the development of quasi-optimal methods for another class of polynomial approximation: non-intrusive interpolation or collocation methods, is an important problem to study. These methods are practical and convenient in that they allow the use of legacy, black-box deterministic numerical solver and the simultaneous approximation of parameterized solutions can be considered as a modular post-processing step. With observation that the accuracy of the interpolation operator $\mathcal{I}_{\Lambda_M}$ is dictated by the inequality 
\begin{align*}
\left\| u - \mathcal{I}_{\Lambda_M}[u]\right\|_{ L^{\infty}(\Gamma)} & \le (1+\mathbb{L}_{\Lambda_M}) \inf_{v\in \mathbb{P}_{ \Lambda_M}} \left\| u - v \right\|_{ L^{\infty}(\Gamma)} 
\\
& \le (1+\mathbb{L}_{\Lambda_M})  \left\| u - u_{\Lambda_M} \right\|_{ L^{\infty}(\Gamma)}, 
\end{align*}
where $\mathbb{L}_{\Lambda_M}$ denotes the Lebesgue constant, we expect that the interpolation schemes in the quasi-optimal subspaces recover the convergence rates described in this work. However, to construct a non-intrusive hierarchical interpolant, two difficult challenges need to be addressed. First, the number of interpolation points needs to remain equal to the dimension of the polynomial space, thus, they must be nested and increase linearly. Second, to guarantee the accuracy of $\mathcal{I}_{\Lambda_M}[u]$, the Lebesgue constant must grow slowly with respect to the total number of collocation points, and we will need to explore the selections of abscissas which optimize this growth.

\begin{acknowledgements}
The authors wish to graciously thank Prof.~Ron DeVore for his interested in our work, his patience in 
discussing the analysis of ''best $M$-term'' approximations, 
and his tremendously helpful insights into the theoretical developments we pursued in this paper.  

This material is based upon work supported in part by the U.S.~Air Force of Scientific Research under grant number 1854-V521-12 and by the U.S.~Department of Energy, Office of Science, Office of Advanced Scientific Computing Research, Applied Mathematics program under contract and award numbers ERKJ259 and ERKJE45; and by the Laboratory Directed Research and Development program at the Oak Ridge National Laboratory, which is operated by UT-Battelle, LLC., for the U.S.~Department of Energy under Contract DE-AC05-00OR22725.
\end{acknowledgements}

\bibliographystyle{spmpsci}      
\bibliography{quasi_optimal_v5_0715}   

\end{document}